\documentclass[preprint,12pt]{elsarticle}
\usepackage{amsmath,amssymb,amsthm,a4wide}
\usepackage[utf8]{inputenc}
\usepackage{amsmath}
\usepackage{subfigure}
\usepackage{amsthm}
\usepackage{amsmath}
\usepackage{amsfonts}
\usepackage{csquotes}
\MakeOuterQuote{"}
\usepackage{tikz} \usepackage{graphicx}
\usepackage{epstopdf}
\usepackage{enumerate}
\usepackage[multiple]{footmisc}
\usepackage{setspace}
\usepackage{amssymb}
\def\@endtheorem{\endtrivlist}%
\newtheorem*{proposition*}{Proposition}
\newtheorem*{question*}{Question}
\newcommand{\phiK}{\phi_{\rm K}}
\newcommand{\phiKbar}{\phi_{\bar{\rm K}}}
\newcommand{\abs}[1] {\left\lvert #1 \right\rvert}

\newcommand{\done}{}

\usepackage{amssymb}

\begin{document}

\begin{frontmatter}
\title{Loss of Physical Reversibility in Reversible Systems}

\author[1,*]{Amir Sagiv}
\author[1]{Adi Ditkowski}
\author[2]{Roy H.\ Goodman}
\author[1]{Gadi Fibich}
\address[1]{School of Mathematical Sciences, Tel Aviv University, Tel Aviv 6997801, Israel}
\address[2]{Department of Mathematical Sciences, New Jersey Institute of Technology, University Heights, Newark, NJ 07102, USA}
\address[*]{corresponding author: asagiv88@gmail.com}

\begin{abstract} 
A dynamical system is said to be reversible if, given an output, the input can always be recovered in a well-posed manner. Nevertheless, we argue that reversible systems that have a time-reversal symmetry, such as the Nonlinear Schr{\"o}dinger equation and the $\phi ^4$ equation can become "physically irreversible". By this, we mean that realistically-small experimental errors in measuring the output can lead to dramatic differences between the recovered input and the original one. The loss of reversibility reveals a natural "arrow of time", reminiscent of the thermodynamic one, which is the direction in which the radiation is emitted outward. Our results are relevant to imaging and reversal applications in nonlinear optics.
\end{abstract}

\begin{keyword}
Nonlinear Schrodinger, Nonlinear dynamics, $\phi^4$ equation, Reversibility, Stability
\end{keyword}

\end{frontmatter}

\section{Introduction}

Consider a time-evolution dynamical system with solution operator $u(t)=Q(t)u(0)$. Strictly speaking, a system is said to be reversible if for any $t>0$ and~$u(t)$~there exists a time-reversal operator $Q^{-1}(t)$ such that $u(0)=~Q^{-1}(t)u(t)$. Thus, a system is reversible if it is always possible to recover the input $u(0)$ from the output $u(t)$. A vast body of research is devoted to proving reversibility in various physical and mathematical systems, see e.g., \cite{alifanov2012inverse, fouque2007wave, lamb1998time}.
 
Obviously, not all dynamical systems are reversible. Indeed, when two different initial states $u_{1}(0)$ and $u_{2}(0)$ evolve in finite time $t_{\rm f}$ into the same output state $u(t_{\rm f})=Q(t_{\rm f})u_1(0)=Q(t_{\rm f})u_2(0)$, then it is impossible to determine the input state from the output state. Hence, the operator $Q^{-1}(t_{\rm f})$ does not exist and so the system becomes irreversible.\footnote{For example, let $y'(t) = -3y^{2/3}$. Then for any $c \in \mathbb{R}$,
 $y_c(t)= \left\{
 \begin{array}{ll}
 (c -t)^3 & {\rm if}~t<c\, ,\\
 0 & {\rm if}~t\geq c \, 
 \end{array}\right.$
is a solution. Hence, if $y(t_{\rm f})=0$ at some time $t_{\rm f}>0$, then $y(0)$ cannot be uniquely determined.}

Usually, for a system to be referred to as reversible, well-posedness of the time-reversal operator $Q^{-1}$ is also required. Intuitively, well-posedness means that for two "close" output states $u_1(t_{\rm f})$ and $u_2(t_{\rm f})$, the corresponding inputs $u_1(0)$ and $u_2(0)$ should be "close" as well.\footnote{A canonical example where the reversal operator exists but is ill-posed is $y'(t)=-y$ with $y(0) = y_0$. Since $y(t) = y_0e^{-t}$, then for any $t_{\rm f}>0$ and $y(t_{\rm f})$, there corresponds a unique $y(0)=y(t_{\rm f})e^{t_{\rm f}}$, and so reversibility is possible. Because of the exponential dependence in $t_{\rm f}$, however, the slightest error in $y(t_{\rm f})$ will have a large impact on the recovered value of $y(0)$. The canonical PDE example for this type of irreversibility is the heat equation, which is well-posed forward in time but ill-posed backward in time \cite{beck1985inverse}.} This requirement guarantees that minor errors in the measurement of the output would not lead to large errors in the recovered input. See e.g.,~\cite{bertero1998introduction, fink2004time, stuart2010inverse} for various methods of reversal and recovery of the input state. 

In this study, we argue that systems that are considered reversible under the above definition (existence and well-posedeness of $Q^{-1}$), can nevertheless exhibit irreversibility in a weaker yet physically meaningful way. Consider, for example, the nonlinear Schr{\"o}dinger equation (NLS) $$
i\partial _t \psi (t,{\bf x})  + \Delta \psi  +N(|\psi|)\psi =0  \, , \qquad 
\psi(0,{\bf x}) = \psi _0 ({\bf x}) \in H^1  \, , \qquad {\bf x}\in \mathbb{R}^d ,$$
{\done where $N(|\psi|):\mathbb{R_+}\to \mathbb{R}$ is real-valued}. The NLS has the time-reversal symmetry $
t \to -t $ and $\psi \to \psi ^{\star}$. Hence, given~$\psi(t_{\rm f},\cdot)$, one can recover the original initial condition~$\psi_0 $~by solving the NLS backward to $t=0$. Furthermore, since the NLS solution is well-posed in $H^1$ (so long as it exists) \cite{strauss1990nonlinear}, then by the time-reversal symmetry it is also well-posed backward in time. Therefore, the NLS is reversible in the sense that $Q^{-1}(t)$ exists and is well-posed. 

Nevertheless, we argue that the NLS can become "physically irreversible". By this, we mean that realistically-small experimental errors in measuring the output $\psi(t_{\rm f}, \cdot)$ can lead to dramatic differences between the recovered input and the original one. This "loss of reversibility" is due to the generic process whereby NLS solutions approach a solitary wave (or multiple solitary waves) on compact domains, while emitting radiation to infinity. Because the system is Hamiltonian, these are not attractors in the usual sense, because the convergence to solitary waves is only on compact domains, while dismissing "far away" and low-energy radiation. Nevertheless, these "quasi-attractors" are the cause for the loss of reversibility, since their basins of attraction contain initial conditions which are very different from each other, yet they all evolve into the same solitary wave (while emitting radiation). Hence, the "reversal information" for recovering the initial condition lies in the far-away low-amplitude radiation, rather than in the high-amplitude solitary waves. Moreover, the low-amplitude radiation undergoes diffraction/dispersion. This makes it prone to inaccurate measurements, which in turn may result in loss of reversibility.

NLS loss of reversibility is thus a consequence of (1) its solitary waves being quasi-attractors under forward propagation and (2) the time-reversal symmetry, since this implies that they are also quasi-attractors under backward propagation. Hence, generically, when an inward radiation interacts with a backward propagating solitary wave, it is likely to remain a solitary wave. Therefore, if the initial condition is quite different from the solitary wave, a successful reversal of the output to the original initial condition requires a precise measurement of the output radiation (which is the cumulative result of all the radiation emitted outward throughout the forward propagation). Consequently, under perturbations of the output radiation, the backward solution is quite likely to stay near the quasi-attractor solitary wave, rather than escape to the original initial condition.

The above explanation for loss of reversibility in the NLS suggests that a system which has a time-reversal symmetry may lose it in the case where its solutions converge to a "quasi-attractor" while emitting radiation. Indeed, we demonstrate a similar loss of reversibility in kink-antikink interactions in the~$\phi ^4$~equation~\cite{Campbell:2019}. {\done The full scope of phenomena in which reversibility might be lost is quite broad; This main examples for loss of reversibility in this paper are dispersive, nonlinear, non-integrable Hamiltonian systems. However, we also present an example of loss of reversibility in hyperbolic conservation laws (Sec.~\ref{sec:burgers} and~\cite{fink2004time}), and an example in an integrable dispersive system (Sec.\ \ref{sec:3waves}). Moreover, while all of the above examples are nonlinear, we expect a similar loss of reversibility in reversible linear systems that evolve into localized stable states while emitting radiation to infinity, e.g., the linear Schr{\"o}dinger equation with a localized potential \cite{hall2013quantum}, a notion known as asymptotic completeness \cite{enss1978asymptotic}. Generically, we expect loss of reversibility whenever a "large" family of initial conditions converge to the same quasi-attractor.}

In our NLS simulations, loss of reversibility can be observed even when the $L^2$ error of the reversed output is as small as $1\%$. This may go against the standard physical intuition, since the $L^2$ norm corresponds to the power of the beam. The above attractor+radiation interaction dynamics clarifies this apparent inconsistency: Since most of the $L^2$ norm/power of the output is concentrated in the high-amplitude attractor, a $1\%$ change in the overall power can correspond to a significant change of the radiation field (which contains the reversal information). 

One might argue that using the $L^2$ error of the output to predict NLS reversibility is wrong, since the theory of NLS well-posedness is usually formulated in $H^1$ space. We observe, however, that small $H^1$ errors of the output may also lead to loss of reversibility. Thus, although the well-posedness theory guarantees reversibility for a sufficiently small $H^1$ error of every output $\psi(t_{\rm f},\cdot)$, it does not predict the {\em size} of this $H^1$ environment. More rigorously, $H^1$ continuity only guarantees that for every $\epsilon >0$, $t>0$, and $u_0\in H^1$, there exists $\delta >0$ such that for every $\|u_1 - u_0\|_{H^1} \leq \delta$ then $\|Q(t)u_1 - Q(t)u_0 \|_{H^1} \leq \epsilon$. It says nothing on the dependence of $\delta$ on the perturbation distance $\epsilon$, time $t$, or the location in phase space $u_0$. As we show in this paper, $\delta$ may be exceedingly small, which leads to a loss of reversibility. Moreover, since the radiation in the NLS disperses, the size of the well-posedness environment decreases with propagation, as we indeed observe numerically. Thus, the existence and well-posedness of the reversal operator $Q^{-1}(t)$ are not sufficient for the NLS to be reversible in various physically-meaningful setups.

Because of its time-reversal symmetry, it is commonly thought that forward and backward propagation in the NLS are physically equivalent (unlike e.g., in the heat equation). The loss of reversibility, however, suggests that the forward and backward directions are {\em not equivalent} in terms of stability under perturbations. This is reminiscent of the thermodynamic arrow of time. Recall the well-known experiment where a small bottle of perfume is opened in a large room. Microscopically, particles evolve by deterministic and reversible interactions. Macroscopically, however, the forward dynamics, in which the particles are distributed evenly in the room, is irreversible, due to the second law of thermodynamics. This seeming contradiction is settled since, on the microscopic level, the probability of the backward process (the gas particles spontaneously return to the bottle) is negligible under random perturbations.\footnote{\done The rigorous derivation of the time irreversible (macroscopic) Boltzmann equation from the reversible equations of motion for the (microscopic) particles is still very much an open question \cite{mischler2012kac}.}  In this sense, the generic evolution into a quasi-attractor+radiation process which leads to loss of reversibility in the NLS suggests that the radiation induces an "arrow of time". This arrow of time may, or may not, coincide with actual time. Thus, if we associate physical well-posedness with time moving forward, than time is forward propagating in the direction in which the radiation is emitted outward, and the recovery of the input data by backward propagation becomes less and less likely as the solution propagates forward. This point of view further shows that the $L^2$ and $H^1$ errors of the output are inadequate indicators for reversibility, since these norms do not take into account the transverse direction of the radiation.

Whether the NLS is physically irreversible is of practical importance. Indeed, in the last decade there have been algorithmic and experimental attempts at holography \cite{barsi2009imaging, barsi2013nonlinear, goy2011digital, goy2013imaging, goy2015improving}, phase retrieval \cite{lu2013phase, lu2016enhanced}, and pulse reconstruction \cite{berti2014reversibility, tsang2003reverse} in focusing Kerr media, all of which rely on the time-reversal symmetry of the NLS. This study, therefore, reveals some of the fundamental limitations that any such reversal technique would face, and highlights the importance of accurately capturing the radiation for reversal experiments to succeed.

The paper is organized as follows. In Sec.\ \ref{sec:lor_fuse} we demonstrate how loss of reversibility occurs in the fusion of two solitary waves of the one-dimensional cubic-quintic NLS under several seemingly small perturbations. Analysis of loss of reversibility and its relation to a natural "arrow of time" in the NLS is presented in Sec.~\ref{sec:arrow}. We present two additional examples of loss of reversibility - self focusing filaments in the two-dimensional cubic-quintic NLS (Sec.\ \ref{sec:lor2d}) and kink-antikink collisions in the $\phi ^4$ equation (Sec.\ \ref{sec:phi4}). In Sec.\ \ref{sec:imag} we discuss the implications of our findings to optical applications and experiments. Sec.\ \ref{sec:burgers} concludes with a comparison to loss of reversibility in nonlinear acoustics.

\section{Loss of reversibility}\label{sec:lor_fuse}
Light propagation is generally considered to be a reversible process. Indeed, mathematical models from ray optics and the wave equation to Maxwell's equations are all invariant under the transformation $t\to -t$. This invariance is in sharp contrast to other physical processes such as heat diffusion, which are {\em not} reversible.

The propagation of high-power laser beams and pulses is described by the Nonlinear Schr{\"o}dinger equation (NLS) in $d+1$ dimensions
\begin{equation}\label{eq:nls_gen}
i\frac{\partial}{\partial z} \psi (z,{\bf x})  + \Delta \psi  + N(|\psi|)\psi =0  \, , \qquad 
\psi(0,{\bf x}) = \psi _0 ({\bf x}) 
\end{equation}
where ${\bf x} = (x_1,\ldots ,x_d)$ are the transverse coordinates (and/or time in the anomalous dispersive regime), $\Delta =~\partial ^2 _{x_1} +~\cdots +~\partial ^2 _{x_d}$, and $z$ is the dimensionless propagation distance. Since mathematically speaking, $z$ is the evolution variable of~ \eqref{eq:nls_gen}, it will be henceforth referred to as "time". When the medium is absorption-free, then $N(|\psi|)$ is real and the NLS \eqref{eq:nls_gen} has the reversal symmetry
\begin{equation}\label{eq:phase_inv}
z \to -z \, , \qquad \psi \to \psi ^{\star} \, .
\end{equation}
In other words, if $\psi(z,{\bf x} )$ is a solution of \eqref{eq:nls_gen}, then so is $\psi^{\star}(-z,{\bf x})$. Intuitively, this reversibility means that propagation in the positive and negative $z$ directions are physically equivalent. Therefore, in principle, given~$\psi(z_{\rm f},{\bf x})$ at some~$z_{\rm f} >0$, one can recover the original initial condition~$\psi_0 ({\bf x})$~by "time reversal" (or "phase conjugation"), i.e., by solving the NLS~\eqref{eq:nls_gen}~{\em backward} until $z=0$.

\begin{figure}[h!]
\centering
	{\includegraphics[scale=.55]{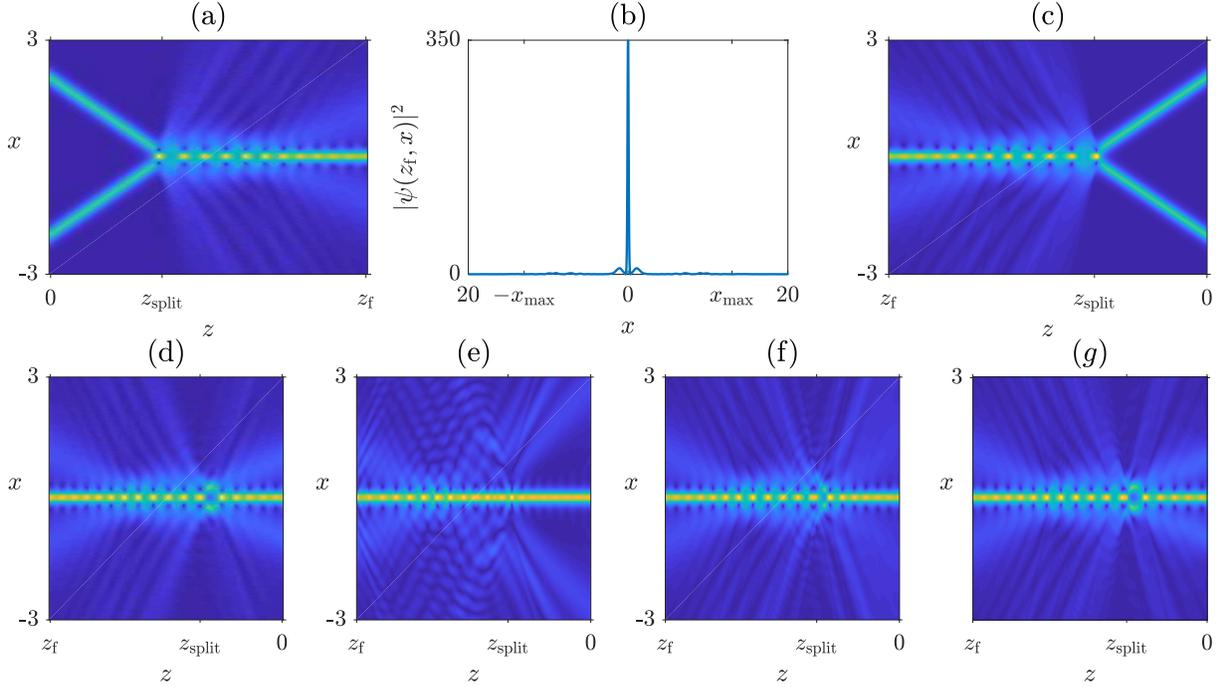}}
\caption{Intensity $|\psi|^2$ of the solution of the one-dimensional cubic-quintic NLS~\eqref{eq:cq_nls} with $\epsilon = 10^{-3}$. (a) Forward propagation of the initial condition \eqref{eq:fusion_ic} for $0\leq z \leq z_{\rm f}$. Here $z_{\rm f} = 0.95$ and $z_{\rm split} = 0.35$. (b) Output intensity $|\psi (z_{\rm f} ,x)|^2$. Here $x_{\max} =13$. (c)~Backward propagation of $\psi(z_{\rm f},x)$ for $z_{\rm f} \geq z \geq 0$. (d--g)~Same for $\psi ^{\rm per} (z_{\rm f},x)$ given by \eqref{eq:truncation}, \eqref{eq:psikmax}, \eqref{eq:phase_mismatch}, and~\eqref{eq:block}, respectively.}
\label{fig:inter_reverse}
\end{figure}

As our first example, consider the one-dimensional cubic-quintic NLS
\begin{equation}
\label{eq:cq_nls}
i\psi_z (z,x)  + \psi_{xx} + |\psi|^2\psi  - \epsilon |\psi|^4\psi =0  \, , \qquad \psi(0,x)=\psi_0 (x) 
\end{equation}
with $\epsilon = 10^{-3}$. This equation admits the solitary wave solutions $ \psi  = e^{i\kappa z} R_{\kappa } (x)  $, where $R_{\kappa}$ satisfies $ -\kappa R_{\kappa } (x) +R''_{\kappa} +R^3 _{\kappa} - \epsilon R_{\kappa } ^5 = 0$. Consider the solution of \eqref{eq:cq_nls} with
\begin{equation}\label{eq:fusion_ic}
\psi _0 (x) = e^{-i\theta x} R_{\kappa} (x-x_0) + e^{i\theta x}R_{\kappa}(x+x_0) \, ,
\end{equation}
where $\kappa = 90$, $x_0=2$, and $\theta = \frac{7}{8}\pi $. This initial condition consists of two intersecting in-phase solitary waves. Upon colliding, the two beams fuse into a single on-axis beam, see Fig.\ \ref{fig:inter_reverse}(a).

If we reverse the "output" beam $\psi (z_{\rm f},x)$, will it "know" that it should split into two solitary waves? Since the NLS \eqref{eq:cq_nls} is reversible, the answer should be positive. Indeed, when we solve \eqref{eq:cq_nls} backward from $z=z_{\rm f}$ to $z=0$, the reversed beam splits into two separate beams, see Fig.\ \ref{fig:inter_reverse}(c).

In physical settings, however, one cannot capture the output beam $\psi (z_{\rm f},x)$ exactly. Hence, one effectively reverses a {\em perturbed} output beam $\psi ^{\rm per} (z_{\rm f},x)$. Two examples of such perturbed profiles are:

\begin{enumerate}
\item {\bf Spatial truncation.} The output beam is measured with a detector of finite radius~$x_{\max}$:\footnote{Reversal can be optically accomplished by phase-conjugators (OPC), e.g., using four-wave mixing in $\chi ^{(3)}$ materials \cite{yariv1977opc}. The same limitations and perturbations may apply to such devices, and we shall therefore refer to them as "detectors" as well.}

\begin{equation}\label{eq:truncation}
\psi ^{\rm per} (z_{\rm f},x) = \left\{ \begin{array}{ll} 
\psi(z_{\rm f},x) \, , &  |x|<x_{\max} \, ,\\ 0 \, , & |x|\geq x_{\max} \, .
\end{array} \, \right. 
\end{equation}

\item {\bf Band-limited filter}. The detector can only resolve spatial frequencies within a band-limited range:
\begin{equation}\label{eq:psikmax}
\psi ^{\rm per} = \mathcal{F}^{-1} \left[ \hat{\psi} ^{\rm BL} (z_{\rm f},k) \right] \, , \qquad  \hat{\psi} ^{\rm BL} (z_{\rm f},k) :\, = \left\{ \begin{array}{ll}
\hat{\psi}(z_{\rm f} ,k) \, , & |k|\leq k_{\max} \, ,\\ 0 \, , & |k|>k_{\max} \, ,
\end{array} \right. 
\end{equation} 
where $\mathcal{F}[\psi (z_{\rm f} ,x)] = \hat{\psi}(z_{\rm f},k)$ is the spatial Fourier transform. This can be the result of a limited resolution (e.g., in CCD cameras), or of a finite-size detector that is placed at a distance after nonlinear media.
\end{enumerate} 

If, instead of reversing the exact output beam $\psi (z_{\rm f},x)$, we reverse the perturbed output $\psi^{\rm per} (z_{\rm f} ,x)$, which is given either by \eqref{eq:truncation} with $x_{\max} =13$ or by \eqref{eq:psikmax} with $k_{\max} = 1.2\pi$, the backwards dynamics changes from beam splitting (Fig.\ \ref{fig:inter_reverse}(c)) to a single on-axis beam (Figs.\ \ref{fig:inter_reverse}(d) and \ref{fig:inter_reverse}(e), respectively). Hence, reversibility is "completely" lost under these perturbations.

\subsection{Indicators for loss of reversibility}

The intensity in the truncated region $|x|\geq x_{\max}$ of the perturbed output beam \eqref{eq:truncation} is {\em $1000$ times smaller} than its peak intensity, see Fig.~\ref{fig:inter_reverse}(b). Hence, it is surprising that truncating this low-intensity region causes such a dramatic loss of reversibility. Intuitively, truncation of a low-intensity region should be "justified" if its power is small compared with the overall power, i.e., if
\begin{equation}
\Delta P : \, = \frac{\|\psi(z_{\rm f} ,\cdot) - \psi ^{\rm per}(z_{\rm f},\cdot)\|_2 ^2}{\|\psi(z_{\rm f} ,\cdot)\|_2 ^2} \ll 1 \, .
\end{equation} 
For both perturbations \eqref{eq:truncation}--\eqref{eq:psikmax}, however, $\Delta P$ is small, as $\Delta P  = 1.8\%$ and $5\%$, respectively. 

The seeming contradiction between the smallness of $\Delta P$ and the complete loss of reversibility might be resolved if we recall that most of the rigorous analytic theory of existence, blowup, and stability in the NLS is carried out in $H^1$ spaces, \cite{strauss1990nonlinear, sulem2007nonlinear, tao2007nonlinear,  fibich2015nonlinear},\footnote{Well-posedness results for {\em some} NLS models \eqref{eq:nls_gen} do exist in $L^2$ spaces \cite{dodson2016global} and in $H^s$ spaces with $s<1$~\cite{killip2010energy}. These theories, however, are not as comprehensive as the $H^1$ theory, and to the best of our knowledge do not exist for \eqref{eq:cq_nls} or for its two-dimensional counterpart \eqref{eq:2d_cq}.} i.e., for solutions with a finite $H^1$ norm, where$$\|\psi\|_{H^1} ^2 : \, = \|\psi\|_2 ^2 + \|\nabla \psi\|_2 ^2  \, .$$ 
Intuitively, the $H^1$ norm is more informative than the power/$L^2$ norm, since it is also affected by the beam's {\em phase}, whereas the power is only affected by the beam's amplitude. Indeed, if we denote $\psi =Ae^{iS}$ where $A$ and $S$ are real, then $$\|\psi \|_{H^1} ^2 - \| \psi  \|_2 ^2 = \|\nabla \psi \|_2 ^2 = \|\nabla A\|_2 ^2 + \|A\nabla S\|_2 ^2 \, .$$ Hence, even in regions where $A$ is moderately small, there can be non-negligible contributions to the $ H^1$ norm from	 $\nabla S$. It is thus more informative to consider the $H^1$ counterpart of $\Delta P$, which we define as\footnote{To maintain consistency with the definition of~$\Delta P$,~we consider here the $H^1$ norm {\em squared}.}

\begin{equation}
\Delta  H^1   : \, = \frac{\|\psi(z_{\rm f} ,\cdot) - \psi ^{\rm per}(z_{\rm f}, \cdot )\|_{H^1} ^2}{\|\psi(z_{\rm f} ,\cdot)\|_{H^1} ^2} \, .
\end{equation}
{\done $\Delta H^1$ is unsuitable for some of the perturbations we consider, e.g., \eqref{eq:truncation}, since $\psi ^{\rm per}$ has a "jump" and is therefore not in $H^1$. However, repeating the simulation in Fig.\ \ref{fig:inter_reverse}(d)  for a variation of \eqref{eq:truncation} where the truncation is replaced by a piecewise linear filter yields loss of reversibility as well (results not shown). We therefore use $\tilde{\Delta}H^1=\|\psi(z_{\rm f},\cdot)\|_{H^1(|x|> x_{\max})}^2 / \|\psi(z_{\rm f},\cdot)\|_{H^1}^2$ for \eqref{eq:truncation}.}
%, and use Parseval identity to define $\tilde{\Delta}H^1 = [\int _{|k|> k_{\max}}(1+|k|^2)|\hat{\psi}(z_{\rm f},k)|^2 \, dk] / \|\psi(z_{\rm f}, \cdot ) \|_{H^1}^2$ for \eqref{eq:psikmax}}

For perturbation, \eqref{eq:truncation} we have that $\tilde{\Delta} H^1 =7.4\%$ (which is still quite small)  whereas for perturbation \eqref{eq:psikmax} we have $\Delta H^1 =29.1\%$. Thus, in both cases $\Delta H^1$ is an order of magnitude larger than $\Delta P$.\footnote{To understand why $\Delta P \ll \Delta H^1$, we note that by Parseval's identity
$\|f \|_2 = \int_{\mathbb{R}} |\hat{f}(k)|^2 \, dk $ and $\|\nabla f\|_2 ^2 = \int_{\mathbb{R}} |k|^2|\hat{f}(k) |^2 \, dk$, where $\hat{f}(k)$ is the Fourier transform of $f(x)$. Hence, the band-limited filter~\eqref{eq:psikmax}~"leaves out" more
$\Delta H^1$ norm then $\Delta P$, since $\int_{k_{\max}}^{\infty} |k|^2 |\hat{f}(k)|^2 \, dk \geq k_{\max} ^2 \int_{k_{\max}}^{\infty} |\hat{f}(k)|^2 \, dk$. The same also applies for the spatial truncation \eqref{eq:truncation}. Since high wave-numbers disperse faster to $x\gg 1$, since far from $x=0$, the intensity is weak, i.e., $|\psi|^2 \ll 1$, then its dynamics are described by the linear Schr{\"o}dinger equation. Hence, the spatial truncation also amounts to the attenuation of high wave numbers.}
Nevertheless, $\Delta H^1$ is still small for perturbation~\eqref{eq:truncation}, yet it leads to loss of reversibility. Hence, even though $\Delta H^1$ is a better indicator for loss of reversibility than $\Delta P$, it is far from providing a definite answer.

{\done Finally, we suggest a possible method to {\em experimentally measure} $\Delta H^1$. Since $\Delta P$ can be measured by intensity ($|\psi|^2$) measurement only, and since phase measurements are difficult in many settings, it would be desirable to suggest an intensity-only measurement method for $\Delta H^1$ as well. First, note that by Parseval Identity $\|f\|_{H^1}^2 = \int _{\mathbb{R}^d} (1+|\xi|^2)|\hat{f}(\xi)|^2 \, d\xi $, where $\hat{f}$ is the Fourier transform of $f$ \cite{hecht2002optics}. Second, note that the linear far-field expression for the propagation of an input profile $f(x)$ is $\hat{f}$. Hence, if the experimental settings allow for {\em linear} propagation of $\psi (z_{\rm f},x)$ far beyond $z_{\rm f}$, then a simple intensity measurement would provide $\|f\|_{H^1} ^2$, as well as $\Delta H^1$.}

\subsection{Digital measurements}

\begin{figure}[h!]
\centering
{\includegraphics[scale=.55]{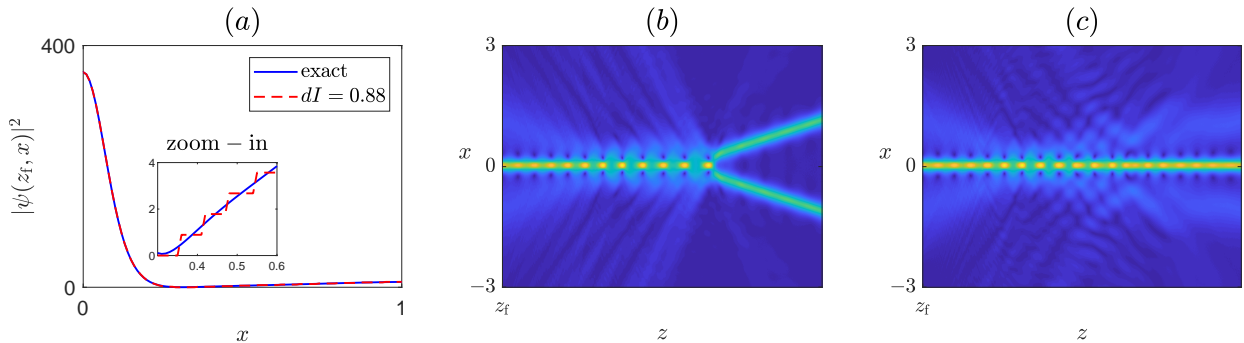}}
\caption{Same settings as in Fig.\ \ref{fig:inter_reverse}. (a) Intensity (solid) and discretized intensity with~$dI=~0.88$ (dash-dots) at $z_{\rm f}$. (b) Backward propagation of $\psi  (z_{\rm f},x)$ discretized with $dI=0.25$. (c)~Same with $dI=0.88$.}
\label{fig:discrete}
\end{figure}

Standard digital detectors use a finite set of discrete values to record the intensity $|\psi |^2$. To model the effect of digitization on back-propagation, we project the output intensity~$|\psi(z_{\rm f},x)|^2$ on the discrete values $I_n = n\cdot dI$ with $dI>0$, and reverse the discretized profile. The discretization with $dI=0.88$ has a seemingly negligible effect on the output profile, see Fig.~\ref{fig:discrete}(a). Indeed, $\Delta P$ for the discretized profiles with $dI=0.88$ and $dI=0.25$ is $4.5\%$ and $1.3\%$, respectively. These two resolutions, however, exhibit very different backward dynamics: The finely-discretized profile undergoes splitting (Fig.\ \ref{fig:discrete}(b)), whereas the coarsely-discretized profile fails to split (Fig.\ \ref{fig:discrete}(c)). Note that in contrast to our prior perturbations,~$\Delta H^1$~is undefined for a digitized output beam.\footnote{The discretized $\psi (z_{\rm f},x)$ is a linear combination of step functions, which are not in $H^1$.}

\subsection{Loss of reversibility and physical ill-posedness}

\begin{figure}[h!]
\centering
{\includegraphics[scale=.4]{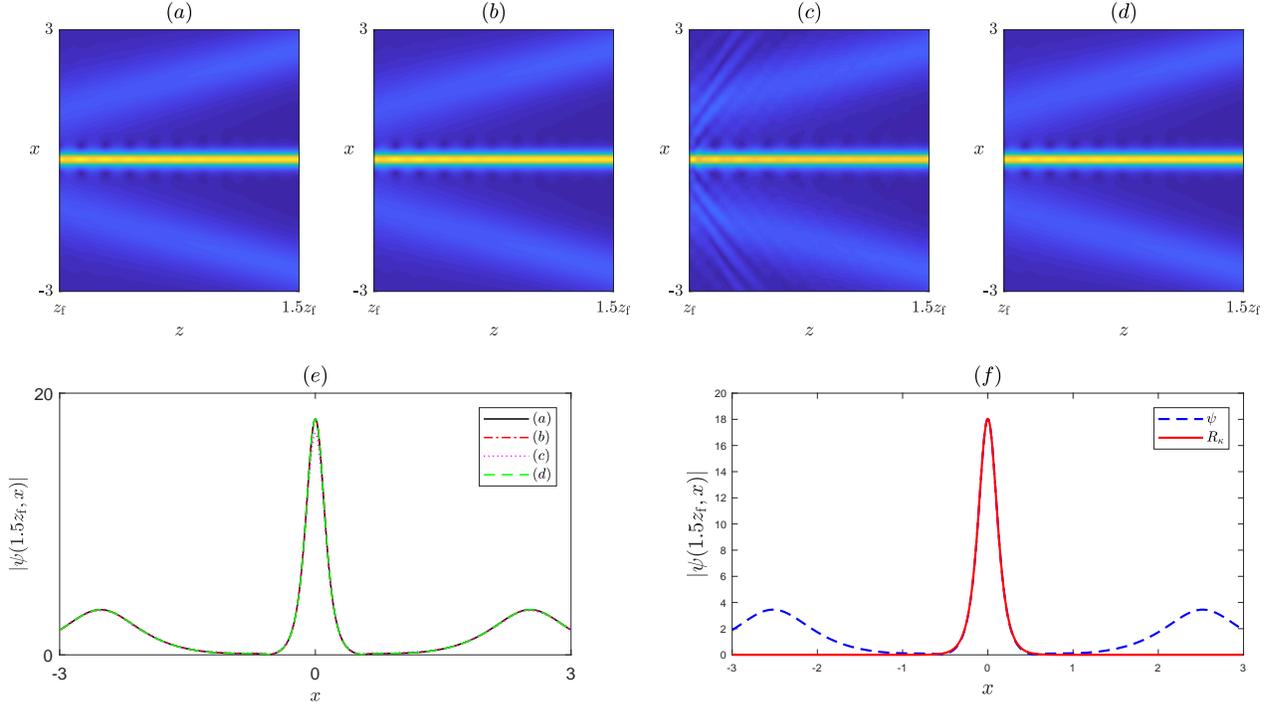}}
\caption{Same settings as in Fig.\ \ref{fig:inter_reverse}. (a). Forward propagation of $\psi(z_{\rm f},x)$ for $z_{\rm f}\leq z \leq~1.5z_{\rm f}$. (b)--(d)~Same for $\psi^{\rm per} (z_{\rm f},x)$ as in Figs.\ \ref{fig:inter_reverse}(d)--\ref{fig:inter_reverse}(f), respectively. (e) $|\psi(1.5 z_{\rm f},x)|$ for subplots~(a)--(d). The four lines are nearly indistinguishable. (f) $|\psi(1.5z_{\rm f},x)|$ (dashes) and the solitary wave $R_{\kappa}(x)$ with $\kappa \approx 127.5$ (solid). The two lines are indistinguishable for~$|x|\leq 1$.}\label{fig:fusion_prop}
\end{figure}

In light of the reversal transformation \eqref{eq:phase_inv}, loss of reversibility is equivalent to physical ill-posedness for the initial condition $\psi ^{\star} (z_{\rm f}, x)$. Here, by {\em physical ill-posedness} we mean that the backward-propagating solution is highly sensitive to small perturbations of its initial condition, i.e., that there exist output profiles $\psi ^{\rm per} (z_{\rm f},x)$ for which $\Delta H^1$ is small, but such that when back-propagated to $z=0$ yield solutions $\psi^{\rm per} (z=0,x)$ that are very different from $\psi_0 (x)$.

The $H^1$ norm is invariant under conjugation, i.e., $\|\psi\|_{H^1} = \|\psi ^* \|_{H^1}$. Therefore, if $\Delta H^1$ is a good indicator for loss of reversibility, then loss of reversibility of $\psi (z_{\rm f}, x)$ should also imply physical ill-posedness in forward propagation of $\psi (z_{\rm f},x)$. To test this hypothesis, in Figs.~\ref{fig:fusion_prop}(a)~-~\ref{fig:fusion_prop}(c) we solve the NLS forward in $z$ for $z\geq z_{\rm f}$ with three initial conditions at $z= z_{\rm f}$: The exact solution $\psi(z_{\rm f},x)$, and the two perturbed profiles $\psi ^{\rm per} (z_{\rm f},x)$ given by~\eqref{eq:truncation}~and~\eqref{eq:psikmax}, respectively. The three solutions are nearly the same at $z=1.5z_{\rm f}$ (Fig. \ref{fig:fusion_prop}(e)), showing that the effect of perturbations \eqref{eq:truncation}~and~\eqref{eq:psikmax} is negligible at $z>z_{\rm f}$. Hence, we again see that $\Delta H^1$ is not a reliable indicator for loss of reversibility. 

\section{Arrow of "time"}\label{sec:arrow}

Figs.\ \ref{fig:inter_reverse} and \ref{fig:fusion_prop} show that the NLS solution with the initial condition $\psi (z_{\rm f},x)$ is physically ill-posed in backward propagation but {\em well-posed} in forward propagation. These seemingly opposite behaviors have a common explanation:

\begin{enumerate}
\item  For $z\geq z_{\rm f}$ (i.e., forward propagation), the high-intensity "core" of $\psi$ is approximately a solitary wave, see~Fig.~\ref{fig:fusion_prop}(f).\footnote{\done Solitary waves of \eqref{eq:2d_cq} are given by $2\kappa ^{1/2}[1+\sqrt{1-\frac{8}{3}\kappa \epsilon}\cosh (2\kappa^{1/2}x)]^{-1/2}  $, see \cite{pushkarov1979self}.} By orbital stability \cite{cazenave1982orbital, weinstein1983nonlinear},~a solution of \eqref{eq:cq_nls} which is close in $H^1$ to a solitary wave will remain so as it propagates. Thus, orbital stability explains the observed well-posedness in forward propagation.\footnote{Fig.\ \ref{fig:fusion_prop}(e) demonstrates the stability of the solution's amplitude. By orbital stability, the {\em complex} profile of the perturbed solution may differ from that of the solitary wave by a constant phase term~$e^{i\beta (z)}$~\cite{cazenave1982orbital, weinstein1983nonlinear}.}

\item For $z\leq z_{\rm f}$ (i.e., backward propagation), since the NLS is continuous with respect to the initial condition, any sufficiently small perturbation of $\psi (z_{\rm f},x)$ would preserve reversibility, and so the back-propagating beam would split. As the perturbation increases, the perturbed profile $\psi ^{\rm per} (z_{\rm f},x)$ enters the "large" $H^1$-neighborhood of the orbit $\{e^{i\beta}R_{ \kappa}\,|\,\beta \in~[0,2\pi] \}$, in which, by orbital stability, the solution remains close to $R_{\kappa}$ for all $z$. In such a case, the back-propagating beam would not split into two beams, and so reversibility would be lost. 
\end{enumerate}

\begin{figure}[h!]
\centering
{\includegraphics[scale=.4]{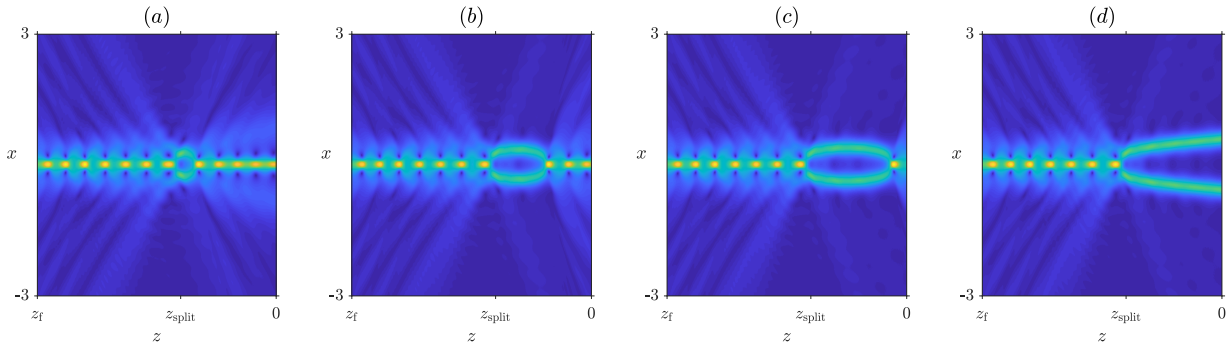}}
\caption{Same settings as in Fig.\ \ref{fig:inter_reverse}(d) with (a) $x_{\max} = 13.6$, (b) $x_{\max} = 13.86$, (c)~$x_{\max}=~13.91$, and (d) $x_{\max} = 14.1$.}
\label{fig:ellipse}
\end{figure}

\begin{figure}[h!]
\centering
{\includegraphics[scale=.5]{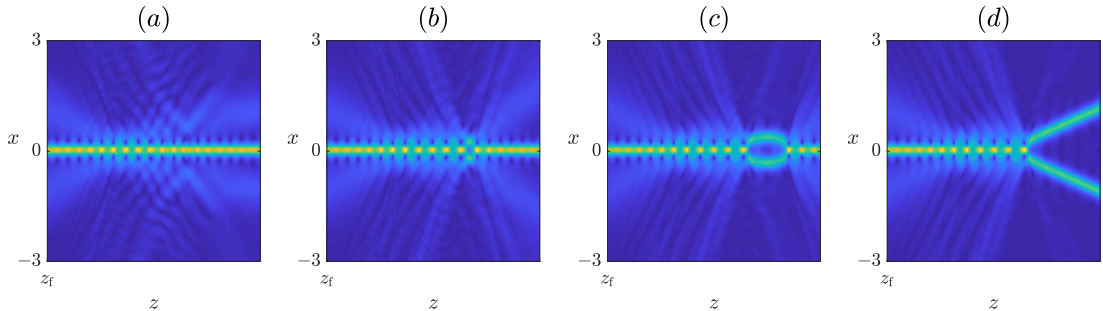}}
\caption{Same settings as in Fig.\ \ref{fig:discrete} with (a) $dI=0.88$, (b) $dI=0.44$, (c)~$dI=0.30$, and (d) $dI=0.24$.}
\label{fig:discrete_ellipse}
\end{figure}

Therefore, in the backward propagation of $\psi(z_{\rm f},x)$, one expects a {\em phase transition} between reversibility (beam splitting) and loss of reversibility (an on-axis solitary wave). To see that, we reconsider the time-reversal of $\psi (z_{\rm f},x)$ under the spatial-truncation \eqref{eq:truncation}, only this time with $x_{\max} = 13.6$ instead of $x_{\max}=13$. At $z=z_{\rm split}$ where the exact solution splits, the perturbed solution develops a double-peak profile, seen in Fig.\ \ref{fig:ellipse}(a) as an ellipse-shaped pattern in the $(z,x)$ plane. As $x_{\max}$ increases to $13.86$ and $13.91$ the ellipse extends further along the $z$ axis, see Figs.\ \ref{fig:ellipse}(b)--\ref{fig:ellipse}(c), respectively. For $x_{\max}= 14.1$, the beam splits and so reversibility is maintained, see Fig.\ \ref{fig:ellipse}(d). Finally, we note that a similar phase transition occurs for other perturbations that lead to loss of reversibility, e.g., as one changes the discretization parameter in (Fig.\ \ref{fig:discrete_ellipse}).

In Figs.\ \ref{fig:inter_reverse} and \ref{fig:fusion_prop} we saw that $\psi(z_{\rm f},x)$ is physically ill-posed in backward propagation but well-posed in forward propagation. Since the NLS is invariant under transformation \eqref{eq:phase_inv}, {\em the opposite} holds for $\psi ^{\star} (z_{\rm f},x)$: It is physically ill-posed in forward propagation but well-posed in backward propagation. Therefore:
\begin{enumerate}
\item NLS well-posedness is not "time-symmetric".

\item The question of ill- and well-posedness is not simply a matter of moving forward or backward in "time". 
\end{enumerate}  

To understand why there is well-posedness in one direction but physical ill-posedness in the opposite direction, we revisit Fig.\ \ref{fig:inter_reverse}(a), where two beams collide and fuse into a single solitary wave. During the fusion process, some radiation is emitted outward. For $z>z_{\rm f}$, the radiation propagates away from the beam, and so perturbations of the surrounding radiation have a minor effect on the dynamics of the beam core. When $\psi (z_{\rm f},x)$ is reversed, it is the radiation that propagates inward that "splits" the solitary wave. From a physical perspective, this observation leads to the following definition of the "arrow of time" - {\em the well-posed forward/ well-posedness direction in the NLS is the one in which the radiation is emitted outward.}

Having an inward-propagating radiation is thus a {\em necessary} condition for the back-propagating beam to split. It is not, however, a sufficient condition. For example, when the same inward-propagating radiation has a $\pi$ phase shift, i.e.,
\begin{equation}\label{eq:phase_mismatch}
\psi ^{\rm per} (x) = \left\{ \begin{array}{ll} 
~~~ \, \psi(z_{\rm f},x) \, , & |x|<x_{\max} \, ,\\ e^{i\pi} \psi(z_{\rm f},x) \, , & |x|\geq x_{\max} \, ,
\end{array} \, \right. 
\end{equation} 
the reversed beam does not split {\done ($\tilde{\Delta} H^1=29\%$ and $\Delta P = 7.2\%$)}, see Fig.\ \ref{fig:inter_reverse}(f).\footnote{\done Here, $\tilde{\Delta}H^1= |1-\beta|\|\psi(z_{\rm f},\cdot ) \|_{H^1(|x|>x_{\max})}^2 / \|\psi(z_{\rm f},\cdot) \|_{H^1}^2$.} Reversibility can also be lost without a phase mismatch between the incoming radiation and the beam core. To see that, we set
\begin{equation}\label{eq:amplify}
\psi ^{\rm per} (x) = \left\{ \begin{array}{ll} 
~\,\psi(z_{\rm f},x) \, , & |x|<x_{\max} \, ,\\ \beta \psi(z_{\rm f},x) \, , & |x|\geq x_{\max} \, ,
\end{array} \, \right.
\end{equation}  
where, as in Fig.\ \ref{fig:inter_reverse}, $x_{\max} = 13$. By continuity, for $\beta$ sufficiently close to $1$ {\done ($0.5 \leq \beta \leq 1.3$) reversibility is not lost. Already for $\beta = 0.4$ ($\tilde{\Delta} H^1 = 2.7\%$ and $\Delta P = 0.6\%$) or $\beta =1.4$ ($\tilde{\Delta} H^1  =1.2\%$ and $\Delta P = 0.3\%$), however,} reversibility is lost and we observe a single backward propagating beam (results not shown).\footnote{\done Here, $\tilde{\Delta}H^1 = \|\psi(z_{\rm f},\cdot) \|_{H^1(|x|\not\in [x_{\rm b}, x_{{\rm b}+1}])}^2/ \|\psi\|_{H^1}^2$.}

We further demonstrate the importance of having the precise radiation with a gentler perturbation of $\psi (z_{\rm f},x)$, in which only a portion of the radiation is blocked, i.e.,
\begin{equation}\label{eq:block}
\psi ^{\rm per} (x) = \left\{ \begin{array}{ll} 
 \psi(z_{\rm f},x) \, , &|x|<x_{\rm b} \, , \\
 0 \, , & x_{\rm b}<|x|<x_{\rm b}+1 \, , \\
 \psi (z_{\rm f},x) \, , &|x|>x_{\rm b} +1 \, .
\end{array} \, \right. 
\end{equation}
For $x_{\rm b}=11$, this is a truly small perturbation ($\Delta P = 0.29\%$, $\Delta H^1 = 3.5\%$). Nevertheless, it prevents beam splitting and thus leads to a complete loss of reversibility, see Fig.\ \ref{fig:inter_reverse}(g). 

The dynamics in Fig.\ \ref{fig:inter_reverse}(c) may give the false impression that the splitting at $z=z_{\rm split}$ in back-propagation is caused by the radiation that "hits" the central beam precisely at $z_{\rm split}$. Applying~\eqref{eq:block}~with $x_{\rm b}=4$,~$8$, and~$10$, however, yields similar dynamics to that of Fig.~\ref{fig:inter_reverse}(g) (results not shown). This shows that "all" of the radiation is needed for the split, and not just the radiation that arrives at~$z_{\rm split}$.

{\bf Remark.} The connection between well-posedness and the direction of the radiation is reminiscent of the Sommerfeld radiation condition in the Helmholtz equation, where solutions are well-posed only if energy does not flow into the system from infinity.

\subsection{Loss of reversibility increases with propagation distance}\label{sec:lor_wtime}

\begin{figure}[h!]
\centering
{\includegraphics[scale=.36]{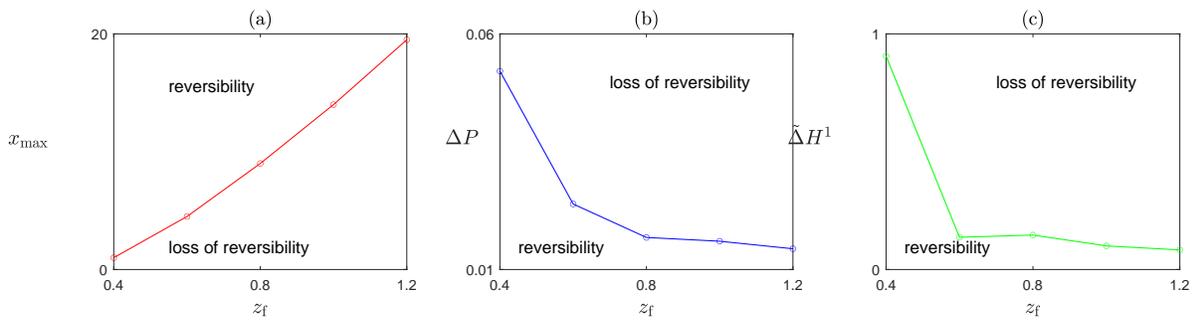}}
\caption{Domains of reversibility and irreversibility for perturbation~\eqref{eq:truncation}. The separatix is (a) $x_{\rm th} (z_{\rm f})$. (b)~$\Delta P(z_{\rm f})$. (c) $\tilde{\Delta} H^1(z_{\rm f})$. }
\label{fig:lor_wtime}
\end{figure} 

So far we saw that for a given one-parameter family of perturbations of $\psi (z_{\rm f},x)$, as the distance between $\psi ^{\rm per} (z_{\rm f},x)$ and $\psi (z_{\rm f},x)$ increases, there is a "phase transition" between reversibility and loss of reversibility. For example, the truncation~\eqref{eq:truncation} at $z_{\rm f}=0.95$ leads to loss of reversibility for $x_{\max} = 13.91$, but not for $x_{\max} =  14.1$ (Fig.\ \ref{fig:ellipse}). {\done For concreteness, we say that reversibility is lost when the output of the reversal is a single beam and not two outgoing beams. We define $x_{\rm th}=x_{\rm th}(z_{\rm f})$ to be the threshold value of $x_{\max}$ such that reversibility is maintained for $x_{\max}>x_{\rm th}$, but is lost for $x_{\max} < x_{\rm th}$.} Fig.\ \ref{fig:lor_wtime}(a) shows that $x_{\rm th}$ increases with $z_{\rm f}$. This result is intuitive, since the low-intensity radiation which contains the reversal information undergoes linear dispersion/diffraction, and so it spreads over a larger spatial domain as $z_{\rm f}$ increases. This argument, however, does not indicate whether the ($L^2$ or $H^1$) distance between $\psi(z_{\rm f},x)$ and $\psi ^{\rm per} (z_{\rm f},x;x_{\max}=~x_{\rm th}(z_{\rm f}))$ increases with $z_{\rm f}$, i.e., whether as $z_{\rm f}$ increases, one is required to capture "more of the radiation" to guarantee reversibility. Figs.\ \ref{fig:lor_wtime}(b) and \ref{fig:lor_wtime}(c) show that both $\Delta P$ and $\Delta H^1$ for perturbation~\eqref{eq:truncation} with $x_{\max}=x_{\rm th}$ decay with $z_{\rm f}$.\footnote{\done Note that the other perturbations employed in Fig. 1 also have a numerical parameter, and one can therefore consider its threshold value in a similar manner.} {\done However, these simulations are inconclusive as to the asymptotic limit, i.e.,:
\begin{question*}
For a given initial condition $\psi _0$, does the minimal $\Delta H^1$ and $\Delta P$ for loss of reversibility of $\psi(z,x)$ vanishes as $z\to \infty$?
\end{question*}
}

\section{Second example - beam collapse}\label{sec:lor2d}

\begin{figure}[h!]
\centering
{\includegraphics[scale=.35]{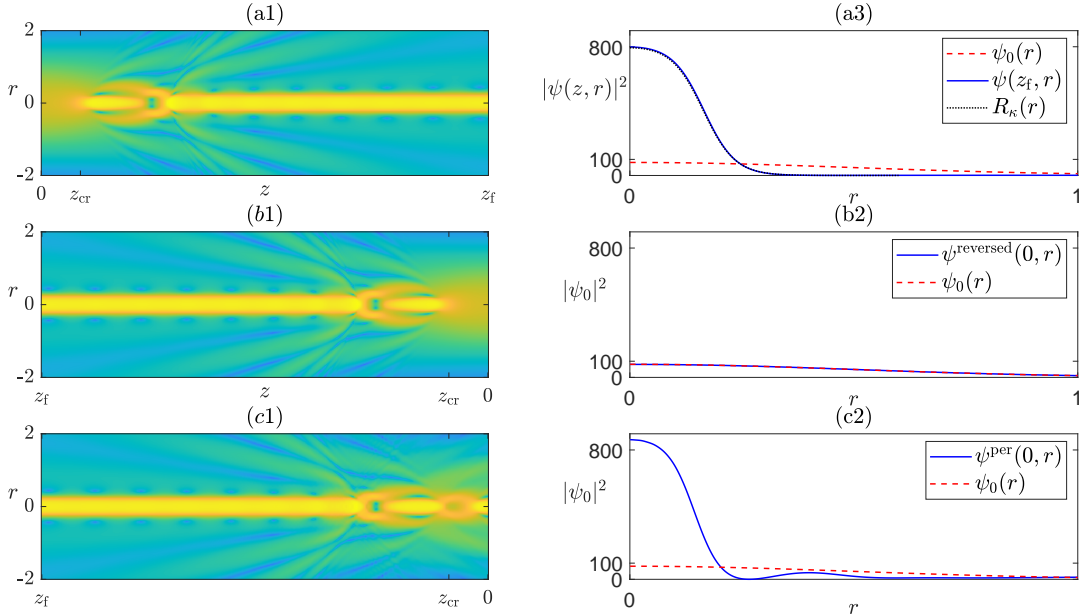}}
\caption{The two-dimensional cubic-quintic NLS \eqref{eq:2d_cq} with the Gaussian input beam \eqref{eq:gaus_ic}. (a1) $|\psi(z,r)|^2$ for $0\leq z \leq z_{\rm f}$. Here $z_{\rm f}=0.4$ and $z_{\rm cr} = 0.04$. (a2) Intensity at $z=0$ (dashes) and at $z_{\rm f}=0.4$ (solid). The dotted line is the solitary wave $R_{\kappa}(r)$ with $\kappa \approx 148$. (b1) Backward propagation of $\psi(z_{\rm f},r)$ for $z_{\rm f}\geq z \geq 0$. (b2) $\psi^{\rm reversed}(z=0,r)$ (solid) is indistinguishable from $\psi _0 (r)$~(dashes). (c1) Same as (b1) for $\psi ^{\rm per} (z_{\rm f},r)$  given by \eqref{eq:truncation} with $r_{\rm max} = 17$. (c2) $\psi ^{\rm per}(0,r)$ (solid), $\psi (z_{\rm f},r)$ (dots), and $\psi _0 (r)$ (dashes).}
\label{fig:gaus_rev}
\end{figure}

Loss of reversibility is not limited to beam fusion. To see that, we solve the two-dimensional cubic-quintic NLS 
\begin{equation}\label{eq:2d_cq}
i\psi_z (z,{\bf x})  + \Delta\psi  + |\psi|^2\psi  - \epsilon |\psi|^4\psi =0  \, , \qquad \psi(0,{\bf x})=\psi_0 ({\bf x})  \, ,
\end{equation} where ${\bf x} =(x,y)$ and $\epsilon = 10^{-3}$, with the single-beam initial condition
\begin{equation}\label{eq:gaus_ic}
\psi_0({\bf x}) = 9e^{-r^2} \, , \qquad r:\,=|{\bf x}| \, ,
\end{equation}
see Fig.\ \ref{fig:gaus_rev}(a). The power of this Gaussian beam is well above the critical power for collapse ($P \approx 7.4 P_{\rm cr}$, see \cite{fibich2000critical}) and therefore it initially collapses at $z_{\rm cr}\approx 0.035$. After the collapse is arrested by the defocusing quintic nonlinearity, the beam evolves into a narrow solitary wave whose peak intensity is $8$ times higher than that of the input beam, see Fig.~\ref{fig:gaus_rev}(a2). We set~$z_{\rm f}=0.4$, which is long after the beam has collapsed.  

The NLS \eqref{eq:2d_cq} is focusing, regardless of whether one moves forward or backward in $z$. Nevertheless, as $\psi(z_{\rm f},r)$ back-propagates it seemingly evolves into the wider initial Gaussian input \eqref{eq:gaus_ic}, see Fig.\ \ref{fig:gaus_rev}(b). As in our previous example, it is the inward-propagating radiation that causes the back-propagating solitary wave to defocus. To see that, we perturb $\psi (z_{\rm f},r)$ with the spatial truncation~\eqref{eq:truncation}~with $r_{\rm max} = 17$,\footnote{$r_{\max}$ plays the same role as $x_{\max}$ in \eqref{eq:truncation}.} and observe that the back-propagated perturbed beam is very different from the initial condition at $z=0$, see Fig.\ \ref{fig:gaus_rev}(c). Thus, this spatial truncation leads to loss of reversibility, even though both $\Delta P = 0.4\%$ and $\tilde{\Delta} H^1 = 7\%$ are small.\footnote{We used a radially-symmetric perturbation in order to show that loss of reversibility is a separate process from symmetry-breaking. Here $\tilde{\Delta}H^1 = \|\psi\|_{H^1(|x|>r_{\max})}^2 / \|\psi\|_{H^1}^2$.}

\subsection{Regaining reversibility}

\begin{figure}[h!]
\centering
{\includegraphics[scale=.55]{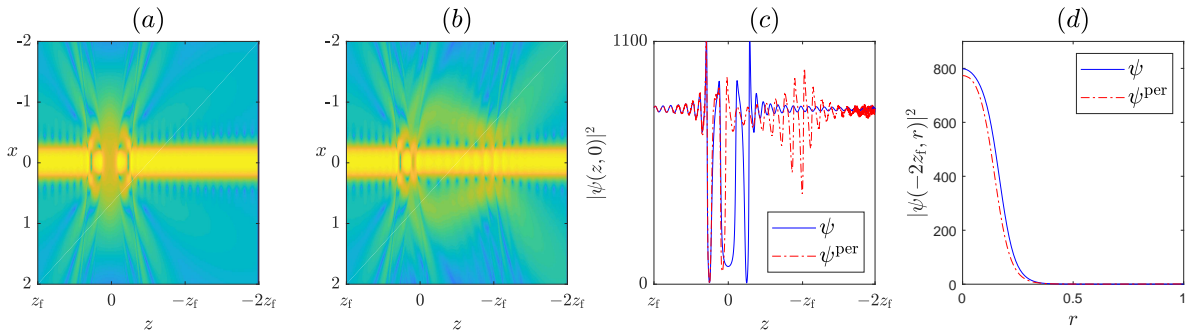}}
\caption{(a)--(b): Same as Figs.\ \ref{fig:gaus_rev}(a) and \ref{fig:gaus_rev}(b), respectively, over a larger domain $z_{\rm f}\geq z \geq-2z_{\rm f}$. (c)~On-axis intensity of the exact (solid) and perturbed (dash-dot) solutions. (d) The output intensity at $z=-2z_{\rm f}$ of the exact (solid) and perturbed (dot dashes) solutions.}
\label{fig:2d_onion}
\end{figure}

In Figs.\ \ref{fig:gaus_rev}(b) and \ref{fig:gaus_rev}(c) we saw that the backward dynamics of $\psi (z_{\rm f} ,x)$ changes dramatically under perturbation \eqref{eq:truncation}, so that the perturbed and unperturbed solutions are very different at~$z=0$. To check whether this loss of reversibility persists beyond~$z=~0$, in Figs.~\ref{fig:2d_onion}(a)~and~\ref{fig:2d_onion}(b) we continue the back-propagation until $z=-2z_{\rm f}$. The two solutions become similar again around~$z~=-1.5~z_{\rm f}$, see Fig.~\ref{fig:2d_onion}(c), so that at $z=-2z_{\rm f}$ they are very close to each other, see~Fig.~\ref{fig:2d_onion}(d).

Why is reversibility regained as $z\to -\infty$? Assume that the perturbation is sufficiently small such that (i) Both the unperturbed and perturbed solutions converge to solitary waves~$R_{\kappa}$~and~$R_{\kappa'}$, respectively, as $z\to-\infty$, and (ii) the perturbed solution emits a similar amount of radiation as the unperturbed solution. Then, since the power of solitary waves of the cubic-quintic NLS \eqref{eq:2d_cq} is monotone in $\kappa$, and since the perturbation is small ($\Delta P \ll 1$), then $|\kappa  - \kappa '| \ll 1$ and so reversibility will be regained as $z\to -\infty$. Note that in this case, as $z\to -\infty$, the "arrow of time" is in the direction of $-z$. Indeed, the radiation is emitted outward for the time-reversal process in this domain.

\section{Third example---kink--antikink collisions in the $\phi^4$ equation}\label{sec:phi4}
Loss of reversibility can also be observed in the $\phi^4$ equation
\begin{equation}\label{eq:phi4}
\phi_{tt}-\phi_{xx} + \phi - \phi^3 = 0 \, .
\end{equation}
In three space dimensions, the $\phi^4$ equation has been used to model physical phenomena at all scales, from the Higgs phenomenon at the subatomic scale~\cite{Higgs:1964} to the formation of domain walls in the early universe at the largest scale~\cite{Zeldovich:1974}, and to phase transitions in superconductors~\cite{Ginzburg:1950}.

Equation \eqref{eq:phi4} possesses a family of traveling-wave kink solutions
\begin{equation*}
%\label{kink}
\phi(x,t)=\phiK(x-v t; v) = \tanh{(\xi/\sqrt{2})} \, , \qquad \xi = (x -x_0- vt)/\sqrt{1-v^2} \, ,
\end{equation*}
for any velocity $-1<v<1$, and another family $\phiKbar:\,=-\phiK$ called antikinks .

\begin{figure}[htbp] %  figure placement: here, top, bottom, or page
   \centering
   \includegraphics[width=.8\textwidth]{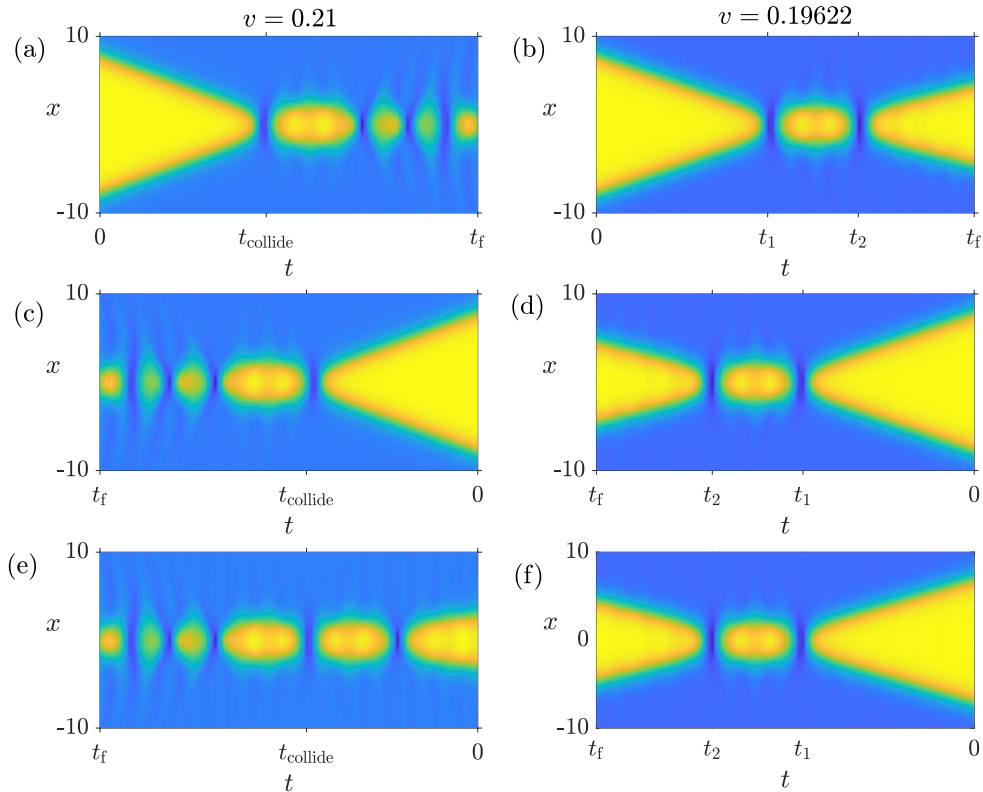} 
   \caption{(a) Numerical simulation of a kink-antikink collision with initial velocity $v=0.21$, showing capture, where $t_{\rm collide}=33$ and $t_{\rm f}=75$. (b) Numerical simulation of kink-antikink collision with $v=0.19622$ showing escape after two collisions. Here, $t_1=34$ and $t_2=52$. (c) Time-reversed simulation from subfigure~(a). (d) Time-reversed simulation from subfigure~(b). (e) Time-reversed simulation from subfigure~(a), data cut off at $x=20$. (f) Time-reversed simulation from subfigure~(b), data cut off at $x=20$. All simulations performed on computational domain~$\abs{x}<32$.}
   \label{fig:phi4}
\end{figure}

Let $\phi_0(x) = \phiK(x+x_0;v) + \phiKbar (x-x_0;-v)+1$  denote a kink $\phiK$ and antikink $\phiKbar$ with equal and opposite velocities $\pm v$; with $\frac{\partial \phi}{\partial t}$ defined analogously. The outcome of their interaction depends very sensitively on the value of $v$ \cite{Campbell:1983em}. For example, for $v=0.21$ the collision results in ``capture'', after which the kink and antikink remain bound together at the collision site and radiate energy away to infinity, see Fig.\ \ref{fig:phi4}(a). For a slightly smaller value of $v$ located inside a specific interval called a "two-bounce window", the kink and antikink first collide, then begin to separate, then reverse direction and collide a second time before escaping with a different final velocity, see Fig.\ \ref{fig:phi4}(d).

We perform a time-reversal experiment, analogous to that described in Sec.~\ref{sec:lor_fuse} for each of the two initial velocities. As expected, running the simulation backward in time, starting at $t_{\rm f}=75$, yields an exact reversal of the initial conditions, see Figs.\ \ref{fig:phi4}(c-d). We then truncate the outputs at $x_{\max}=20$, see \eqref{eq:truncation}, and time-reverse the perturbed outputs, see Figs.~\ref{fig:phi4}(e-f). The truncated ``captured'' solution from Fig.\ \ref{fig:phi4}(a) loses its reversibility, as the recovered input remains captured, rather separating into an escaping kink and anti--kink pair. The truncated two-bounce resonant solution from Fig.\ \ref{fig:phi4}(b), by contrast, maintains its reversibility, since the recovered input is nearly indistinguishable from the original one.\footnote{In this case, reversibility is maintained even if $x_{\max} = 10$ (results not shown).} In both cases, we also simulated the system forward in time and observed that truncating the radiation had a negligible effect on the forward dynamics (results not shown). Thus, similarly to the NLS, the transverse direction of the radiation induces a natural arrow of time with regard to sensitivity to perturbations.

Like NLS, the $\phi ^4$ equation is time reversible, as it is invariant under the transformation $ t\to -t $. It is also well posed for $\left(\phi(\cdot,0),\phi_t(\cdot,0)\right) \in H^2_{\rm loc}\times H^1$, decaying sufficiently rapidly to $\pm1$ as $\abs{x}\to\infty$, see~\cite{Henry:1982kn}. Hence, the reversal operator $Q^{-1}(t)$ exists and is well-posed for all $t>0$. Henry et al. showed that kink solutions are asymptotically orbitally stable, its solutions might become "physically irreversible".

To show that, as in the case of the NLS, the truncated radiation in Figs.\ \ref{fig:phi4}(e) and \ref{fig:phi4}(f) is "small", we need to estimate its relative size. Unlike the NLS, we cannot use the $L^2$ and $H^1$ norms since the solutions of interest do not decay at infinity. However, equation~\eqref{eq:phi4} conserves an energy  $
H :\,= \int_{-\infty}^{\infty} 
(
 \phi_t^2 /2
+ \phi_x^2 /2 
+  \left(\phi^2 - 1 \right)^2/4  \, \, \,dx
)
$,
which is finite for these solutions. We therefore define the relative error as
$$
\Delta H^{\rm rel} :\, = \frac{H(\phi^{\rm per}) - H(\phi)}{H(\phi)}.
$$
The relative errors of the truncated solutions in both cases are small: $0.42\%$ for the "captured" solution and $0.28\%$ for the "two-bounce" solution.

The explanation for the different behaviors is related to the mathematical theory underlying the sensitive dependence of the final state on the initial velocity, as demonstrated in Fig~\ref{fig:phi4}(a-b). We describe this briefly; for further details, see \cite{Goodman:2005vv, Goodman:2019, Sugiyama:1979gl, Takyi:2016eo, Campbell:2019}. The kink and antikink attract each other, and one can define a potential energy describing their interaction, as well as a kinetic energy. While the combined kinetic and potential energy is negative, the kink-antikink pair is bound together. When the combined energy is positive, they escape from each other. During collisions, there are three important effects: \textbf{(1)}~The kink-antikink separation undergoes large acceleration. \textbf{(2)}~The kink and antikink reversibly interchange energy with a secondary mode of oscillation of the system.%
\footnote{The identity of this secondary mode, long thought to be a so-called internal mode, remains an open question. Under this assumption, Sugiyama developed a finite-dimensional model that thoroughly analyzed by Goodman and Haberman.~\cite{Sugiyama:1979gl,Goodman:2005vv}. Takyi and Weigel has shown a major flaw in this reasoning, including an algebra error in Sugiyama's model, rendering its use invalid for quantitative arguments.~\cite{Takyi:2016eo}. Others have suggested it is a quasinormal mode~\cite{Dorey2018}. Nonetheless, as a qualitative description, the model gives excellent insights.}
This mode is unexcited before the first collision, and the amount of energy exchanged on subsequent interactions depends sensitively on its amplitude and phase at the moment of collision. \textbf{(3)}~The creation of radiation, which irreversibly carries energy away from the localized solutions. Energy transferred to the internal mode may be returned to the kinks as kinetic energy, allowing them to escape, but energy lost as radiation cannot, because it carries energy with it away from the kink location toward infinity. The radiation approximately satisfies the linearized evolution, which is dispersive, and high-frequency radiation moves at unit speed, much faster than the kink and antikink.  The time between collision $n$ and $(n+1)$ (and thus the phase of the second oscillator) depends on the combined kinetic and potential energy in the kink-antikink pair following collision $n$. If this energy is positive following a collision, the kink and antikink escape. As more collisions that occur, more energy is lost to radiation and the probability of eventual escape decreases.

\begin{figure}[htbp] %  figure placement: here, top, bottom, or page
  \center
   \includegraphics[width=0.8\textwidth]{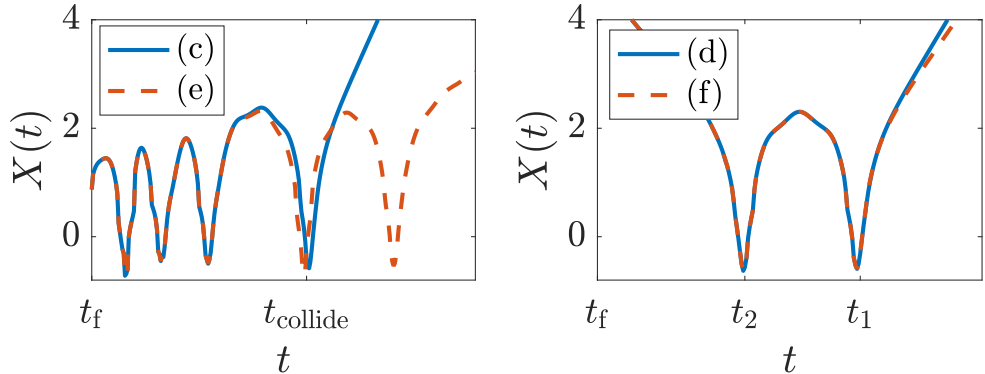} 
   \caption{The location of the antikink over time in the time-reversed simulations shown in Fig.~\ref{fig:phi4}. Panel (a) shows the location for the "capture" simulations and panel (b) for the "two-bounce" simulations.}
   \label{fig:phi4positions}
\end{figure}

Figure~\ref{fig:phi4positions} shows the location of the antikink over time in the four time-reversed runs, and demonstrates how the locations diverge over time in the cutoff and non-cutoff cases.  Just before $t_{\rm collision}$ in the two time-reversed simulations of the "capture" solution, shown in subfigure~(a), the two simulations have begun to diverge. The cutoff simulation has lost a little bit of energy and the collision occurs slightly before the collision in the non-cutoff simulation. Because the result of the collision depends sensitively on the phase of the secondary oscillator, this leads to escape in one simulation but not in the other, thus time reversal is lost. By contrast, the positions in the two reverse-time simulations of the two-bounce solution, shown in subfigure~(b), only begin to diverge after the kink and antikink have already escaped, leading to preservation of time-reversal.

\section{Imaging and reversibility in nonlinear optics}\label{sec:imag}

In recent years, there has been an increasing interest in imaging in nonlinear media. We briefly note four directions of related research.

Berti et al.\ \cite{berti2014reversibility} studied the reversibility of ultra-short pulses in {\em lossy} focusing medium. The envelopes of these pulses, modeled either by the $(3+1)d$ NLS or by more comprehensive models, are numerically shown to be reversible in physical settings and parameters. The authors of \cite{berti2014reversibility} explain the seeming contradiction between reversibility and effects such as energy loss, intensity clamping, loss of phase etc., as follows: {\em "... the information redistribution in space and over all observables ensures that the information required to back-propagate the pulse is in fact conserved. Indeed, it is well-known that the filaments are not an isolated system, but are in strong interaction with the surrounding photon bath"}~\cite{berti2014reversibility}.

Our study supports the conclusion that the information in the surrounding "photon bath" may be crucial for maintaining reversibility. Moreover, we show that imprecise or incomplete "knowledge" of the phase and amplitude of the "photon bath" (i.e., the radiation) may lead to loss of reversibility. Hence, our study suggests that the reversibility which was numerically observed in \cite{berti2014reversibility} might be lost in physical experiments due to e.g., the use of a detector with a finite size or a finite bandwidth. {\done The stabilizing role of the "photon bath" in propagation in nonlinear focusing media has been noticed in several experimental settings, see e.g., \cite{courvoisier2003ultra, kolesic2004self, skupin2004interaction}.
}

Goy and Psaltis \cite{goy2011digital, goy2013imaging} developed an algorithm for holography, (i.e., the recovery of an object's $3D$ structure), in focusing nonlinear media. In their experiment, an object reflects a laser beam which propagates to the detector in nonlinear focusing medium. The imaging algorithm, based on Psaltis' pulse-reversal algorithm \cite{tsang2003reverse} has essentially two steps: (1) measuring the phase and amplitude of the output profile, and (2) numerically solving the NLS backward. At moderate powers, the imaging is improved by using their nonlinear holography process, compared to linear Schr{\"o}dinger-based algorithm. At high powers, however, the success of their algorithm decreases. This failure is attributed to the emergence of "parasitic filaments" in the numerical backward propagation, which obscure the original input image \cite{goy2015improving}. Goy, Makris, and Psaltis proposed an algorithm to improve digital holography in focusing media by introducing random artificial perturbations to the output signal \cite{goy2015improving}. This method successfully addresses the issue of parasitic filaments, and thus improves digital holography in focusing media.

Our study suggests a different limitation to reversibility in focusing media: As the solutions converge to a quasi-attractor (solitary waves), the "reversal information" which is contained in the radiation, disperses. This issue has not been accounted for in \cite{goy2011digital, goy2013imaging}, and so the algorithm in \cite{goy2015improving} is therefore not designed to handle it.

Barsi, Wan, and Fleischer \cite{barsi2009imaging} and Barsi and Fleischer \cite{barsi2013nonlinear} demonstrated that the use of a {\em defocusing} nonlinear medium can improve the imaging resolution beyond Abbe's diffraction limit. Subsequently, the same group harnessed medium nonlinearity to the well-known Gerchberg-Saxton algorithm for phase retrieval, i.e., the retrieval of complex phase using only intensity measurement \cite{lu2013phase, lu2016enhanced}. The authors noted that {\em "... Although focusing nonlinearities can also couple these modes, noise-induced instabilities can dominate the signal and may limit the ability to invert..."} \cite{barsi2009imaging}. Indeed, in the phase-retrieval numerical simulations, focusing media lead to instabilities \cite{lu2016enhanced}. Our work identifies a different inherent limitation to achieving reversibility in focusing media (dispersion of the reversal information contained in the radiation), which may further limit its use for imaging and phase retrieval.

Our work may also be relevant imaging in nonlinear {\em inhomogeneous} media. To see that, we note the work of Frostig et al \cite{frostig2017focusing}. which characterizes the propagation of speckled light(a wide beam with each pixel given an iid random phase) in focusing media by two processes: (i) self focusing. (ii) the formation of many beams and their subsequent fusion into few intense filaments. Combined with our study, these insights suggest that propagation-reversal of light from an inhomogeneous media/speckled source via focusing media will be a daunting task; It will require to reverse both beam fusion and beam collapse - the very two processes that we demonstrated to be prone to loss of reversibility (Sec.\ \ref{sec:lor_fuse} and \ref{sec:lor2d}). 

\subsection{An integrable example: three-wave interaction}\label{sec:3waves}
Loss of reversibility may also occur in beam fusion in {\em three-wave interaction} in quartic nonlinear media ($\chi ^{(2)} \neq 0$), which is described by the {\em integrable} equations \cite{kaup1976three, kaup1979space}
\begin{alignat*}{2}
&\partial _z u_{1}(z,x) + c_1 \partial _x u_1 =&&\gamma _1 u_2 ^* u_3 ^* \, , \\
&\partial _z u_{2}(z,x) + c_2 \partial _x u_2 =  &&\gamma _2 u_3 ^* u_1 ^* \, , \\
&\partial _z u_{3}(z,x) + c_3 \partial _x u_3 =  &&\gamma _3 u_1 ^* u_2 ^* \, ,
\end{alignat*} 
where for $i=1,2,3$, $u_i(z,x)$ is an envelope of a wave with frequency $\omega _i$ and group velocity~$c_i$, and $\gamma _i=\pm 1$ is a medium-dependent coefficient. This is because, on one hand, these equations admit solutions where $u_1$ and $u_2$ solitary waves collide and fuse into a $u_3$ solitary wave while emitting radiation. On the other hand, these equations also admit traveling-waves solutions~$u_3 = f(x-c_3 z)$, where $u_1\equiv u_2 \equiv 0$. Since for both of these solutions, the output at large $z$ consists of a single solitary wave with frequency $\omega_3$ and some radiation, the reversal information is contained in the radiation, and therefore the three-waves interaction is prone to loss of reversibility. 

\subsection{Loss of reversibility in the integrable NLS?}

\begin{figure}[h!]
\centering
{\includegraphics[scale=.55]{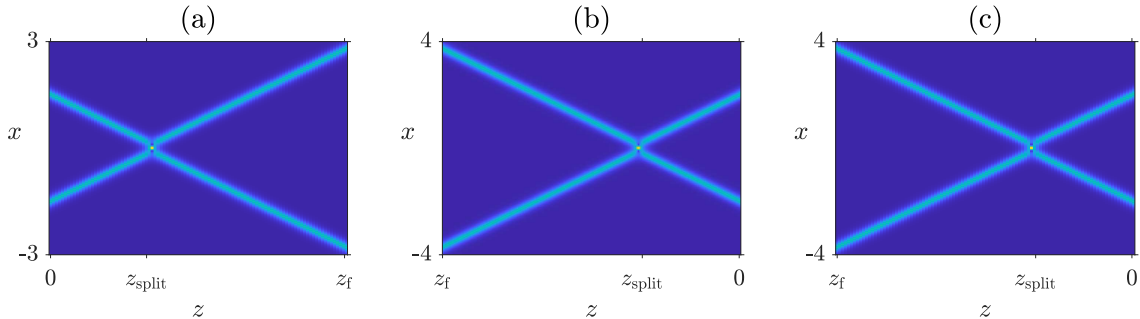}}
\caption{(a)--(b): Same as Figs.\ \ref{fig:gaus_rev}(a) and \ref{fig:gaus_rev}(b), respectively, over a larger domain $z_{\rm f}\geq z \geq-2z_{\rm f}$. (c)~On-axis intensity of the exact (solid) and perturbed (dash-dot) solutions. (d) The output intensity at $z=-2z_{\rm f}$ of the exact (solid) and perturbed (dot dashes) solutions.}
\label{fig:cubic}
\end{figure}
{\done The previous section shows that integrable PDEs can undergo loss of reversibility. Since Sec.\ \ref{sec:lor_fuse} considers the one-dimensional cubic-quintic NLS \eqref{eq:cq_nls}, it is only natural to ask whether the integrable cubic one-dimensional NLS, i.e., \eqref{eq:cq_nls} with $\epsilon= 0$, can undergo a similar loss of reversibility. When we repeat the experiment of propagating, reversing, and reversing a truncated output \eqref{eq:truncation} as in Sec.\ \ref{sec:lor_fuse} but with $\epsilon=0$, the results are substantially different; see Fig.\ \ref{fig:cubic}. Indeed, as predicted by the inverse scattering transform, the number of beams at~$t=\pm \infty$~is conserved, and therefore there is no beam fusion \cite{ablowitz2004ist}. 
Therefore, since we started with a (nearly) pure two-soliton solution, there is a negligible amount of radiation and so the truncation has a negligible effect on reversibility. In the terminology used throughout this paper, there is only one quasi-attractor, the two-soliton solution, and so reversibility is not lost (or even mildly affected) given the truncation. This is not to say that loss of reversibility is impossible in integrable systems (Sec.\ \ref{sec:3waves}) or that it is impossible in the cubic NLS in different settings.}

\section{Comparison with loss of reversibility in Burgers equation}\label{sec:burgers}

It is instructive to compare the loss of reversibility in the NLS and $\phi ^4$ equation with the one in Burgers equation. Loss of reversibility in acoustics was {\em experimentally} demonstrated by Tanter et al.\ \cite{tanter2001breaking}. In their experiments, they observed that reversibility is possible if one reverses the wave before an acoustic shock wave forms, but not if it is measured after the shock.

Tanter et al.\ modeled their experiment with two coupled Burgers equations. To explicitly demonstrate loss of reversibility in the inviscid Burgers equation \begin{equation}\label{eq:burgers} u_t (t,x) + uu_x = 0 \, , \qquad u(0,x)=~u_0(x) \, ,
\end{equation} we consider the two initial conditions
\begin{equation}\label{eq:shock}
u_0 ^{(1)}(x)= \left\{ \begin{array}{ll} 1 \, , & x<0 \, ,\\ 0 \, ,& x\geq 0 \, , \end{array}\right.  \qquad u_0 ^{(2)}(x)= \left\{ \begin{array}{ll} 1 \, ,& x<0 \, ,\\ 
1-x \, , & 0\leq x <1 \, ,\\
0 \, , & x\geq 1 \, . \end{array}\right.   
\end{equation}

\begin{center}
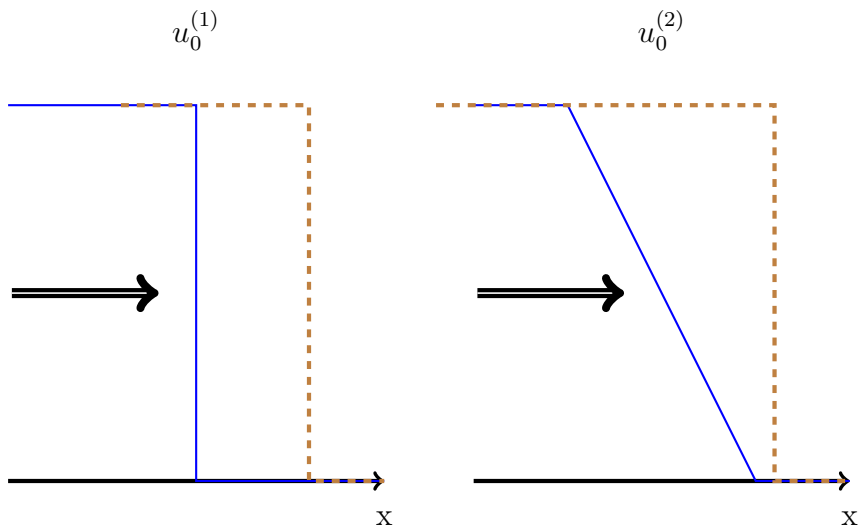
\begin{figure}[h!]
\centering
\subfigure
{
\begin{tikzpicture}[scale=5]
\draw [-> , ultra thick] (-1/2,0) -- (1/2,0);
\draw [blue, thick]
(-1/2,1) -- (0,1) -- (0,0) -- (1/2,0);
\draw [brown, ultra thick, dashed]
(-1/2+.3,1) -- (0+.3,1) -- (0+.3,0) -- (1/2,0);
\draw [double , -> , ultra thick]
(-1/2+.01, 1/2)--(-1/5+.1,1/2);
\node  at (1/2,-.1) {x};
\node at (0,1.2) {$u_0 ^{(1)}$};	
\end{tikzpicture}
}
\subfigure
{
\begin{tikzpicture}[scale=5]
\draw [-> , ultra thick] (-1/2,0) -- (1/2,0);
\draw [blue, thick]
(-1/2,1) -- (-1/4,1) -- (1/4,0) -- (1/2,0);
\draw [brown, ultra thick, dashed]
(-1/2-.1,1) -- (0+.3,1) -- (0+.3,0) -- (1/2,0);
\draw [double , -> , ultra thick]
(-1/2+.01, 1/2)--(-1/5+.1,1/2);
\node  at (1/2,-.1) {x};
\node at (0,1.2) {$u_0 ^{(2)}$};	
\end{tikzpicture}
}
\caption{The inviscid Burgers equation \eqref{eq:burgers}, with two initial conditions as given in \eqref{eq:shock} (solid), and evolves into the same shock at~$t=1$ (dashes) },
\end{figure}
\end{center}
Both initial conditions evolve into the same moving-shock solution $u(t,x)=u_0^{(1)}(x-t)$ for $t\geq 1$ \cite{leveque1990conservation}. Hence, reversibility is lost for $t\geq 1$. Indeed, since a shock wave forms when multiple characteristic lines coincide, its time-reversal (a rarefaction wave) consists of a cone in the $(t,x)$ plane where there is a lacuna of characteristics  \cite{whitham2011linear}. 

The above explicit example of loss of reversibility is different from the previous examples in this paper, since it occurs when two different initial conditions evolve into {\em exactly} the same output profile (rather than into two similar profiles). This is irreversibility in the strictest mathematical sense, since for $t\geq 1$, the inverse operator $Q^{-1}(t)$ does not exist.

\subsection{Open question}
The above "loss of existence" of the reversal operator $Q^{-1}$ is due to the formation of a singularity in Burgers equation (as a shock wave is formed, its derivative becomes infinite). Is there a similar "loss of existence" of $Q^{-1}$ in the NLS? For this to occur, the NLS solution should become singular. Recall that the cubic two-dimensional cubic NLS
\begin{equation}\label{eq:nlsc}
i\psi _z (z,x,y) + \Delta \psi + |\psi |^2 \psi =0 \, ,
\end{equation}
has a "large" family of initial conditions that all lead to blowup solutions which decompose at the singularity point $z_{\rm cr}$ into a universal blowup profile $\psi_{R_0}$ that collapses as a $\delta$ function at the $\log \log$ rate, and a radiation field $\phi$, such that $\|\psi_{R_0}(z,\cdot)\|_{H^1}, \|\psi(z,\cdot)\|_{H^1} \to \infty$ and $\psi-\psi_{R_0} \to \phi$ in $L^2$ as $z\to z_{\rm cr}$ \cite{merle2004universality, fibich2015nonlinear}.
It is currently an open question whether two different initial conditions can collapse at the same point $z_{\rm cr}$ with exactly the same radiation field $\phi$. In other words, is the radiation field $\phi$ at the blowup point $z_{\rm cr}$ sufficient to recover the input beam $\psi _0$? {\done For the special case of $\phi=0$, reversibility can be lost at collapse, as the following proposition shows.
\begin{proposition*}
For every $x_* \in \mathbb{R}^2$ there exists an initial condition $\psi_{0,x_*} \in C^{\infty}(\mathbb{R}^2)$ for which $$\lim\limits_{t\to 1^-} |\psi_{x_*}(t,x)| = \delta(x) \, ,$$
where the limit is in the weak-$*$ topology, and such that the following pointwise convergence holds $$\phi (x) :\,= \lim\limits_{t\to 1^-} \psi_{x_*}(t,x) = 0 \, . \qquad  x\neq 0 \, .$$
\end{proposition*}
See \ref{ap:collapse} for the proof. Unfortunately, this proof give no clue regarding the general case of $\phi\neq 0$, since it relies on the non-generic property of the explicit collapse solution \eqref{eq:exp}; see \cite{fibich2015nonlinear} and the references therein.} This question is also open for other nonlinear PDEs that have blowup solutions which decompose into a universal blowup profile and radiation, such as the nonlinear heat equation \cite{giga1985heat,  fermanian2000stability} or the generalized Kortewege-de Vries equation~(gKdV)~\cite{martel2014blow}.

Physically, singularities never form in acoustics and in nonlinear optics. In the NLS \eqref{eq:2d_cq}, we observed that small collapse-arresting defocusing nonlinearity leads to what we called physical loss of reversibility. Similarly, Tanter et al.\ observed numerically and experimentally that this physical loss of reversibility occurs in the viscous Burgers, where shock formation is arrested~\cite{tanter2001breaking}. One might surmise that the introduction of a dissipative term, which renders the inverse Burgers equation ill-posed, is the source for loss of reversibility. By continuity, however, any sufficiently small regularizing term, dissipative or not, should lead to similar loss of reversibility near the shock formation, i.e., where $Q^{-1}$ ceases to exist in the inviscid case.

\section{Acknowledgments}
The authors would like to thank B.\ Ilan and P.\ D.\ Miller for useful ideas and conversations. The work of G.F.\ and A.S.\ was partially supported by the Israel Science
Foundation (ISF) under grant number 177/13. This research was carried out during a long-term stay of A.S.\ as a guest of Prof.\ R.R.~Coifman and the Department of Mathematics at Yale University, whose hospitality is gratefully acknowledged.

\appendix 
\section{Proof of loss of reversibility for a collapsing beam}\label{ap:collapse}
{\done Denote the explicit blowup solution of the NLS \eqref{eq:nlsc} by
\begin{equation}\label{eq:exp}
\psi_{t_{\rm c}}(t,x) = \frac{1}{|t_{\rm c} -t|} R(\xi)\exp\left\{ i\left( \zeta (t) -\frac{|x|^2}{4(t_{\rm c}-t)}\right)\right\} \, ,  \qquad t_{\rm c} \geq 0 \, ,
\end{equation}
where $R$ is the ground state solution of $-R+\Delta R + R^3 = 0$, and $$\xi = \frac{x}{t_{\rm c} - t} \, , \qquad \zeta = \frac{t}{t_{\rm c}(t_{\rm c} - t)} \, .$$
Fixing $t_{\rm c}=1$, this solution of \eqref{eq:nlsc} satisfies the desired blowup properties \cite{fibich2015nonlinear}.

Next, define the Galilean transformation by
\begin{equation}\label{eq:gal}
\mathcal{F}_c \psi(t,x)  :\,= \psi(t,x-ct)\exp \left\{i\left(\frac{c\cdot x}{2} - \frac{|c|^2z}{4} \right) \right\} \, , \qquad c\in \mathbb{R}^2 \, .
\end{equation}
By the Galilean and Translation Invariances of the NLS \eqref{eq:nlsc} then $\mathcal{F}_{x_*} \psi_{1}(t=0,x-x_*)$ collapses with exactly the same properties at $(t_{\rm c}=1, x = 0)$.}
\bibliographystyle{abbrv}

\begin{thebibliography}{10}

\bibitem{ablowitz2004ist}
M.J.~Ablowitz, B.~Prinari, and A.D.~Trubatch. \newblock {\em Discrete and Continuous Nonlinear Schr{\"o}dinger Systems,} 
\newblock Cambridge University, 2004.

\bibitem{alifanov2012inverse}
O.~Alifanov.
\newblock {\em Inverse {H}eat {T}ransfer {P}roblems}.
\newblock Springer Science \& Business Media, 2012.

\bibitem{Anninos:1991vo} P.~Anninos, S.~Oliveira, and R.~A.~Matzner.
\newblock Fractal structure in the scalar $\lambda  (\varphi ^2-1)^2$ theory.
\newblock {\em Phys. Rev. D,} 44:1147--1160, 1991.




\bibitem{barsi2013nonlinear}
C.~Barsi and J.~W. Fleischer.
\newblock Nonlinear abbe theory.
\newblock {\em Nat. Photon.} 7:639, 2013.

\bibitem{barsi2009imaging}
C.~Barsi, W.~Wan, and J.~W. Fleischer.
\newblock Imaging through nonlinear media using digital holography.
\newblock {\em Nat. Photon.} 3:211, 2009.

\bibitem{beck1985inverse}
J.~Beck, B.~Blackwell, and C.~Clair~Jr.
\newblock {\em Inverse {H}eat {C}onduction: {I}ll-{P}osed {P}roblems}.
\newblock James Beck, 1985.

\bibitem{bertero1998introduction}
M.~Bertero and P.~Boccacci.
\newblock {\em Introduction to {I}nverse {P}roblems in {I}maging}.
\newblock CRC, 1998.


\bibitem{berti2014reversibility}
N.~Berti, W.~Ettoumi, J.~Kasparian, and J.-P. Wolf.
\newblock Reversibility of laser filamentation.
\newblock {\em Opt. Exp.}, 22:21061--21068, 2014.

\bibitem{Campbell:2019}\newblock D.~K.~Campbell.
\newblock {\em Historical Overview of the Model}
\newblock appeared in {\em A Dynamical Perspective on the $\phi ^4$ Model: Past, Present and Future}, edited by P.~G.~Kevrekidis and J.~Cuevas-Maraver, Springer Nature, 2019.

\bibitem{Campbell:1983em} D.~K.~Campbell, J.~F.~Schonfeld, and C.~A.~Wingate.
\newblock Resonance structure in kink-antikink interactions in theory.
\newblock {\em Phys. D,} 9:1-32, 1983.


\bibitem{cazenave1982orbital}
T.~Cazenave and P.-L. Lions.
\newblock Orbital stability of standing waves for some nonlinear
  schr{\"o}dinger equations.
\newblock {\em Comm. Math. Phys.} 85:549--561, 1982.

\bibitem{courvoisier2003ultra}
F.~Courvoisier, V.~Boutou, J.~Kasparian, E.~Salmon,
G.~M{\'e}jean, J.~Yu, and J.P.~Wolf.
\newblock Ultraintense light filaments transmitted through clouds.
\newblock {\em Appl. Phys. Lett.} 83:213--215,~2003.

\bibitem{dodson2016global}
B.~Dodson.
\newblock Global well-posedness and scattering for the defocusing, l
  2-critical, nonlinear {S}chr{\"o}dinger equation when $d=1$.
\newblock {\em Amer. J. of Math.} 138:531--569, 2016.

\bibitem{Dorey2018}
P.~Dorey and T.~Romanczukiewicz.
\newblock Resonant kink--antikink scattering through quasinormal modes. 
\newblock {\em Phys. Lett. B,} 779:117--123, 2018.

\bibitem{enss1978asymptotic}
V.~Enns.
\newblock Asymptotic completeness for quantum mechanical potential scattering.
\newblock {\em Comm. Math. Phys.} 61:285--291, 1978.

\bibitem{fermanian2000stability}
C.~Fermanian~Kammerer, F.~Merle, and H.~Zaag.
\newblock Stability of the blow-up profile of non-linear heat equations from
  the dynamical system point of view.
\newblock {\em Mathematische Annalen}, 317:347--387, 2000.

\bibitem{fibich2015nonlinear}
G.~Fibich.
\newblock {\em The {N}onlinear {S}chr{\"o}dinger {E}quation}.
\newblock Springer, New York, 2015.

\bibitem{fibich2000critical}
G.~Fibich and A.~L. Gaeta.
\newblock Critical power for self-focusing in bulk media and in hollow
  waveguides.
\newblock {\em Opt. Lett.} 25:335--337, 2000.

\bibitem{fink2004time}
M.~Fink, G.~Montaldo, and M.~Tanter.
\newblock Time reversal acoustics.
\newblock In {\em IEEE {U}ltrasonics {S}ymposium}, volume~2, pages 850--859.
  IEEE, 2004.

\bibitem{fouque2007wave}
J.~Fouque, J.~Garnier, G.~Papanicolaou, and K.~Solna.
\newblock {\em Wave {P}ropagation and {T}ime {R}eversal in {R}andomly {L}ayered
  {M}edia}, volume~56.
\newblock Springer Science \& Business Media,~2007.

\bibitem{frostig2017focusing}
H.~Frostig, E.~Small, A.~Daniel, P.~Oulevey, S.~Derevyanko, and Y.~Silberberg.
\newblock Focusing light by wavefront shaping through disorder and
  nonlinearity.
\newblock {\em Optica}, 4:1073--1079,~2017.

\bibitem{giga1985heat}
Y.~Giga and R.~Kohn.
\newblock Asymptotically self-similar blow-up of semilinear heat equations.
\newblock {\em Comm. Pure Appl. Math.} 38:297--319, 1985.

\bibitem{Ginzburg:1950}V.~L.~Ginzburg and L.~D.~Landau.
\newblock {\em Zh. Eksp. Teor. Fiz.} 20:1064, 1950.

\bibitem{Goodman:2019}\newblock R.~H.~Goodman.
\newblock {\em 
Mathematical Analysis of Fractal Kink-Antikink Collisions in the Model.}
\newblock appeared in {\em A Dynamical Perspective on the $\phi ^4$ Model: Past, Present and Future}, edited by P.~G.~Kevrekidis and J.~Cuevas-Maraver, Springer Nature, 2019.

\bibitem{Goodman:2005vv} R.~H.~Goodman, and R.~Haberman.
\newblock Kink--antikink collisions in the equation: The $n$-bounce resonance and the separatrix map.
\newblock {\em SIAM J. Appl. Dyn. Sys.} 4:1195--1228, 2005.


\bibitem{goy2011digital}
A.~Goy and D.~Psaltis.
\newblock Digital reverse propagation in focusing Kerr media.
\newblock {\em Phys. Rev. A}, 83:031802, 2011.

\bibitem{goy2013imaging}
A.~Goy and D.~Psaltis.
\newblock Imaging in focusing Kerr media using reverse propagation.
\newblock {\em Photon. Res.} 1:96--101, 2013.

\bibitem{hall2013quantum}
B.C.~Hall.
\newblock {\em Quantum Theory for Mathematicians.}
\newblock Springer, New York, 2013.

\bibitem{hecht2002optics}
E.~Hecht.
\newblock {\em Optics,}
\newblock Addison-Wesley, Reading, MA, 2002.

\bibitem{Henry:1982kn}
D.~B. Henry, J.~F. Perez, and W.~F. Wreszinski.
\newblock {Stability theory for solitary-wave solutions of scalar field
  equations}.
\newblock {\em Commun. Math. Phys.}, 85:351--361, 1982.

\bibitem{Higgs:1964}
P.~W.~Higgs.
\newblock Broken symmetries and the masses of gauge bosons.
\newblock {\em Phys. Rev. Lett.} 13:508--509, 1964.

%\bibitem{karpman1981pert}
%V.I.~Karpman and V. V.~Solov'eV. 
%\newblock A perturbational approach to the two-soliton systems,
%\newblock {\em  Phys. D,} Nonlinear Phenomena 3:487--502, 1981.


\bibitem{kaup1976three}
D.~Kaup.
\newblock The three-wave interactionâ a nondispersive phenomenon.
\newblock {\em Stud. Appl. Math.} 55:9--44, 1976.

\bibitem{kaup1979space}
D.~Kaup, A.~Reiman, and A.~Bers.
\newblock Space-time evolution of nonlinear three-wave interactions. i.
  interaction in a homogeneous medium.
\newblock {\em Rev. Mod. Phys.} 51:275, 1979.


\bibitem{killip2010energy}
R.~Killip and M.~Visan.
\newblock Energy-supercritical NLS: Critical $H^s$-bounds imply scattering.
\newblock {\em Comm. Part. Diff. Eq.} 35:945--987, 2010.

\bibitem{kolesic2004self}
M.~Kolesic and J.M.~Moloney.
\newblock Self-healing femtosecond light filaments.
\newblock {\em Opt. Lett.} 29:590--592, 2004.

\bibitem{lamb1998time}
J.~Lamb and J.~Roberts.
\newblock Time-reversal symmetry in dynamical systems: a survey.
\newblock {\em Phys. D}, 112:1--39, 1998.

\bibitem{leveque1990conservation}
R.~LeVeque.
\newblock {\em Numerical {M}ethods for {C}onservation {L}aws}.
\newblock ETH Z{\"u}rich, Switzerland, 1990.

\bibitem{lu2013phase}
C.-H. Lu, C.~Barsi, M.~O. Williams, J.~N. Kutz, and J.~W. Fleischer.
\newblock Phase retrieval using nonlinear diversity.
\newblock {\em Appl. Opt.} 52:D92--D96, 2013.

\bibitem{lu2016enhanced}
J.~Lu, C.~Li, and J.~Fleischer.
\newblock Enhanced phase retrieval using nonlinear dynamics.
\newblock {\em Opt. Exp.} 24:25091--25102, 2016.

\bibitem{martel2014blow}
Y.~Martel, F.~Merle, and P.~Rapha{\"e}l.
\newblock Blow up for the critical generalized Korteweg--de Vries equation. i:
  Dynamics near the soliton.
\newblock {\em Act. Math.} 212:59--140, 2014.

\bibitem{merle2004universality}
F.~Merle and P.~Raphael.
\newblock On universality of blow-up profile for l2 critical nonlinear
  Schr{\"o}dinger equation.
\newblock {\em Inven. Math.} 156:565--672, 2004.

\bibitem{mischler2012kac}
S.~Mischler and C.~Mouhot.
\newblock Kac's program in kinetic theory.
\newblock {\em Inven. Math.} 193:1--147, 2012.

\bibitem{goy2015improving}
D.~Psaltis. A.S.~Goy, K.G.~Makris.
\newblock Improving the quality of filament-impaired images in {K}err media by
  statistical averaging.
\newblock {\em Opt. Exp.} 23:431--444, 2015.

\bibitem{pushkarov1979self}
K.~Pushkarov, D.~Pushkarov, and I.~Tomov.
\newblock Self-action of light beams in nonlinear media: soliton solutions.
\newblock {\em Opti. Quant. Elec.} 11:471--478, 1979.

\bibitem{skupin2004interaction}
S.~Skupin, L.~Berg{\'e}, U.~Peschel, and F.~Lederer.
\newblock Interaction of femtosecond light filaments with obscurants in aerosols.
\newblock {\em Phys. Rev. Lett.} 93:023901, 2004.

\bibitem{strauss1990nonlinear}
W.~Strauss.
\newblock {\em Nonlinear {W}ave {E}equations}.
\newblock American Mathematical Soc., 1990.

\bibitem{stuart2010inverse}
A.~Stuart.
\newblock Inverse {P}roblems: a {B}ayesian {P}erspective.
\newblock {\em Act. Num.} 19:451--559, 2010.

\bibitem{Sugiyama:1979gl} T.~Sugiyama.
\newblock Kink--antikink collisions in the two-dimensional model.
\newblock {\em Prog. Theor. Phys.} 61:1550--1563, 1979.


\bibitem{sulem2007nonlinear}
C.~Sulem and P.~Sulem.
\newblock {\em The {N}onlinear {S}chr{\"o}dinger {E}quation: {S}elf-{F}ocusing
  and {W}ave {C}collapse}.
\newblock Springer, 1999.

\bibitem{Takyi:2016eo} I.~Takyi, and H.~Weigel.
\newblock Collective coordinates in one-dimensional soliton models revisited.
\newblock {\em Phys. Rev. D}, 94:085008, 2016.

\bibitem{tanter2001breaking}
M.~Tanter, J.~Thomas, F.~Coulouvrat, and M.~Fink.
\newblock Breaking of time reversal invariance in nonlinear acoustics.
\newblock {\em Phys. Rev. E}, 64:016602, 2001.

\bibitem{tao2007nonlinear}
T.~Tao, M.~Visan, and X.~Zhang.
\newblock The nonlinear {S}chr{\"o}dinger equation with combined power-type
  nonlinearities.
\newblock {\em Comm. Part. Diff. Eq.} 32:1281--1343, 2007.

\bibitem{tsang2003reverse}
M.~Tsang, D.~Psaltis, and F.~G. Omenetto.
\newblock Reverse propagation of femtosecond pulses in optical fibers.
\newblock {\em Opt. Lett.} 28:1873-1875, 2003.

\bibitem{weinstein1983nonlinear}
M.~Weinstein.
\newblock Nonlinear {S}chr{\"o}dinger equations and sharp interpolation
  estimates.
\newblock {\em Comm. Math. Phys.} 87:567--576, 1983.

\bibitem{whitham2011linear}
G.~B. Whitham.
\newblock {\em Linear and Nonlinear Waves}, volume~42.
\newblock Wiley \& Sons, 2011.

\bibitem{yariv1977opc}
A.~Yariv and D.~Pepper.
\newblock Amplified reflection, phase conjugation, and oscillation in
  degenerate four-wave mixing.
\newblock {\em Opt. Lett.} 1:16--18, 1977.

\bibitem{Zeldovich:1974} Y.~B.~Zel'dovich, I.~Y.~Kobzarev, and L.~B.~Okun.
\newblock Cosmological consequences of spontaneous violation of discrete symmetry.
\newblock {\em Zh. Eksp. Teor. Fiz.} 40:3--11, 1974.

\end{thebibliography}

\pdfinclusioncopyfonts=1

\end{document}